\documentclass[a4paper]{article}

\usepackage{booktabs}
\usepackage[english]{babel}
\usepackage[utf8]{inputenc}
\usepackage{amssymb}
\usepackage[fleqn]{amsmath}
\usepackage{arydshln}
\usepackage{amsthm}
\usepackage{graphicx}
\usepackage{framed}
\usepackage[colorinlistoftodos]{todonotes}
\usepackage{subcaption}
\usepackage{multirow}
\usepackage{enumitem}
\usepackage{titling} 
\usepackage{titlesec}
\usepackage[top=1.25in, bottom=1.25in, left=1.25in, right=1.25in]{geometry}
\usepackage{breqn}
\usepackage[numbers]{natbib}
\usepackage{tikz}
\usepackage{pgfplots}
\usepackage{ulem}
\usepackage{caption}
\usetikzlibrary{shapes,snakes}
\usepackage{stackengine}
\usetikzlibrary{arrows,%
	petri,%
	topaths}%
\usepackage{algorithm}
\usepackage{algpseudocode}
\usepackage{tabstackengine}

\TABbinary
\newcounter{example}
\newenvironment{example}[2]{
\refstepcounter{example}
\label{#1}\begin{leftbar}\textbf{\noindent \hspace{-5.3mm}Example \theexample} \\#2\end{leftbar}
}{}

\newcounter{drafts}

\newcounter{module}
\newcounter{model}
\newcommand \linedabstractkw[2]{
  \renewcommand\maketitlehookd{%
    \mbox{}\medskip\par
    \centering
    \hrule\medskip
    \begin{minipage}{0.9\textwidth}
    \textbf{Abstract}\\ #1\\
    
    \textit{Keywords: }#2
    \end{minipage}\medskip\hrule\medskip
    }      
}

\newcommand \ERPWtemplate[3]{
\usepackage{fancyhdr}
\pagestyle{fancy}
\fancyhf{}
\fancyhead[RO]{#1}
\fancyhead[LO]{#2}
\fancyfoot[CO]{#3}
\fancyfoot[RO]{\thepage}

\renewcommand{\headrulewidth}{1pt}
\renewcommand{\footrulewidth}{1pt}
}

\stackMath

\newcommand{\re}[1]{\text{Re}\left(#1\right)}

\newcommand{\w}{\omega}
\newcommand{\ww}{\overline{\omega}}
\newcommand{\f}{\varphi}
\newcommand{\ff}{\overline{\varphi}}

\newcommand{\T}{\mathbb{T}}

\DeclareMathOperator{\tr}{tr}
\DeclareMathOperator{\ND}{ND}

\DeclareMathOperator{\IG}{IG}

\theoremstyle{plain}
\newtheorem{theorem}{Theorem}[section]

\theoremstyle{plain}
\newtheorem{lemma}[theorem]{Lemma}

\theoremstyle{plain}
\newtheorem{proposition}[theorem]{Proposition}

\theoremstyle{remark}

\theoremstyle{plain}
\newtheorem{corollary}{Corollary}[theorem]

\theoremstyle{definition}
\newtheorem{definition}{Definition}

\theoremstyle{plain}

\DeclareRobustCommand{\VAN}[3]{#2} 

\ERPWtemplate{P. Wissing \& E.R. van Dam}{Gain graphs with two eigenvalues}{\leftmark}

\tikzset{group/.style = {shape=circle,draw,dotted,minimum size=1em}}
\tikzset{vertex/.style = {shape=circle,draw,minimum size=1em}}
\tikzset{arc/.style = {->,> = latex'}}
\tikzset{edge/.style = {-,> = latex'}}
\tikzset{negarc/.style = {->,> = latex',dashed}}
\tikzset{negedge/.style = {-,> = latex',dashed}}
\tikzset{tree/.style = {-,> = latex',line width=.7mm}}

\title{
Unit gain graphs with two distinct eigenvalues and systems of lines in complex space
}
\author{
Pepijn Wissing\thanks{Corresponding author: p.wissing@tilburguniversity.edu},~ Edwin R. van Dam \\ \small{Department of Econometrics and Operations Research, Tilburg University}
}

\begin{document}
\linedabstractkw{
Since the introduction of the Hermitian adjacency matrix for digraphs, interest in so-called complex unit gain graphs has surged. 
In this work, we consider gain graphs whose spectra contain the minimum number of two distinct eigenvalues. 
Analogously to graphs with few distinct eigenvalues, a great deal of structural symmetry is required for a gain graph to attain this minimum. 
This allows us to draw a surprising parallel to well-studied systems of lines in complex space, through a natural correspondence to unit-norm tight frames.  
We offer a full classification of two-eigenvalue gain graphs with degree at most $4$, or with multiplicity at most $3$.
Intermediate results include an extensive review of various relevant concepts related to lines in complex space, including SIC-POVMs, MUBs and geometries such as the Coxeter-Todd lattice, and many examples obtained as induced subgraphs by employing a technique parallel to the dismantling of association schemes. 
Finally, we touch on an innovative application of simulated annealing to find examples by computer. 
}{
Complex unit gain graphs, Spectrum, Hermitian, Equal-Norm Frames
}
\maketitle

\section{Introduction}                       
\label{sec: intro}
An especially captivating line of research considers graphs, whose associated matrices have few distinct eigenvalues. 
Such graphs are generally highly structurally symmetric, which allows for a beautiful interplay of algebra and combinatorics. 
A non-empty graph must always have at least two distinct eigenvalues; a bound that is essentially only attained in a complete graph. 
In this work, we explore the degree to which this holds true for modern alternatives to the classical binary graphs, and investigate the necessary circumstances for such generalizations to yield examples outside of their immediate graph parallels. 

In recent times,  algebraic graph theorists have been actively considering various generalizations of the traditional graphs.
The most general such form, known as \textit{complex unit gain graphs} \cite{zaslavsky1989biased, reff2012spectral}, considers what is essentially a weighted bidirected graph with arc weights (known as \textit{gains}) from the complex unit circle, such that the gain of every arc is the complex conjugate of the gain of its converse arc. 
The corresponding \textit{gain matrix}, constructed in the natural way, is clearly Hermitian, and therefore has real eigenvalues. 
Some interesting, relevant works in this currently particularly active field include \cite{belardo2020balancedness, lu2021complex, xu2020complex}. 

A number of special cases of unit gain graphs have been researched to a varying degree by various authors. 
Mostly pioneered  by Zaslavsky, the notion of \textit{signed graphs}, which in hindsight are simply \textit{real} unit gain graphs, has been the topic of ongoing inquiries for several decades. 
See \cite{zaslavsky2012mathematical} for a complete bibliography, and \cite{belardo2019open} for an interesting collection of open problems.
More recently, \citet{guo2017hermitian} and \citet{liu2015hermitian} independently constructed a Hermitian adjacency matrix for mixed graphs.
Its definition is such that every edge is assigned gain $1$ and every (outgoing) arc is assigned gain $i:=\sqrt{-1}$. 
By characterizing mixed graphs in this way, several often-used tools (such as eigenvalue interlacing and the quotient matrix) became available for the directed graph paradigm; see \cite{wissing2019negative, gavrilyuk2019multiplicities} for some related work.
A new variant to this Hermitian adjacency matrix, which employs the sixth, rather than the fourth roots of unity, was also considered by \citet{mohar2020new}.
Building on these ideas, \citet{wissing2020spectral} launched an investigation into \textit{signed directed graphs}, obtained by equipping the edges and arcs of a mixed graph with a sign function.
A natural Hermitian matrix to represent such graphs is then constructed analogously to \cite{guo2017hermitian, mohar2020new} by associating every sixth root of unity with one of the ways in which a vertex can be incident to a (signed) edge or arc.

Recently, Belardo et al. \cite{belardo2019open} posed the problem of investigating signed graphs with exactly two distinct eigenvalues.
Quite a few papers have since appeared on the topic.
In particular, \citet{huang2019induced} has used a construction of signed $n$-cubes with exactly two eigenvalues in his recent proof of the so-called Sensitivity Conjecture of Nisan and Szegedy on Boolean functions.
Furthermore, \citet{ramezani2020some} applies the star-complement technique to find infinitely many $k$-regular signed graphs with two distinct eigenvalues $\pm \sqrt{k}$, with $k=5,6,\ldots,10$, and \citet{stanic2020spectra} offers various theoretical and computational results, among others completing the list of $3$- and $4$-regular signed graphs with two distinct eigenvalues. 
Lastly, in an earlier work classifying cyclotomic matrices, \citet{greaves2012cyclotomic} has obtained several infinite such families with two eigenvalues, slightly  restricted versions of which can be interpreted as unit gain graphs.

In this work, we will further develop the ideas and results on signed graphs, above, to the more general setting of complex unit gain graphs.
Many (or in fact all) such graphs that have exactly two eigenvalues correspond to interesting systems of lines in complex space, such that the angle between every non-orthogonal pair of lines is equal to some given constant. 
If every such line is represented by a vector with a given norm (say, 1), then one obtains an object known as an equal-norm tight frame \cite{waldron2018introduction}.
Moreover, if no two such vectors are orthogonal, the system is said to be equiangular. 

Due to their rich theoretical properties and their numerous practical applications, equiangular tight frames are arguably the most important class of finite-dimensional tight frames, and they are the natural choice when one tries to combine the advantages of orthonormal bases with the concept of redundancy provided by frames \cite{strohmer2003grassmannian}.
While most research regarding equiangular lines is relatively old, the quantum computing community has been increasingly interested in equiangular tight frames, especially in the context of \textit{symmetric, informationally complete, positive operator-valued measures} (SIC-POVM), which is a prominent candidate for a "standard quantum measurement." 
Such SIC-POVM's are equivalent to equiangular tight frames of $d^2$ vectors in $\mathbb{C}^d$, and their existence for arbitrary $d$ is one of the important open problems of the moment in quantum computing.

Our ultimate goal is to classify various families of unit gain graphs, with two distinct eigenvalues. 
The applied approach is twofold.
Specifically, the combinatorially oriented graph perspective is focused on the degree of said graphs, while the lines perspective, that focuses on the multiplicities of eigenvalues, is more algebraically oriented. 
For gain graphs of degree at most four, we are able to completely characterize the collection of desired unit gain graphs. 
Some of these collections have infinitely many switching-distinct members, for given order and degree. 
The lines perspective also produces an abundance of interesting examples, and a complete characterization with least multiplicity at most $3$ is obtained. 
Moreover, various other examples stemming from well-known combinatorial objects such as the Coxeter-Todd lattice are discussed, as well as a technique that is parallel to the dismantling of association schemes, which is used to find many two-eigenvalue subgraphs.

This paper is organized as follows. 
In Section \ref{sec: prelim}, we provide a thorough introduction of the concepts used throughout. 
Section \ref{sec: constructions} is concerned with the recursive constructions that may be applied to grow arbitrarily large gain graphs with exactly two distinct eigenvalues. 
Next, Section \ref{sec: systems of lines} draws from the literature on systems of lines in complex space to construct various examples, and offer the necessary insight to classify all two-eigenvalues gain graphs with small multiplicity. 
Section \ref{sec: small degree} provides classifications of unit gain graphs with restricted degree, taking the combinatorial perspective. 
Finally, in Section \ref{sec: heuristic}, we briefly touch on an application of simulated annealing to search for the desired graphs by computer.

\section{Preliminaries}                      
\label{sec: prelim}
Let us first thoroughly define the key concepts and notation that is used throughout this work.

\subsection{Basic definitions}
\label{subsec: basic def}
Let $G=(V,E)$ be a bidirected graph, whose vertex set $V$ is of order $n$, and whose arc set $E$ consists of ordered pairs of vertices called \textit{arcs}, which are denoted $uv$, for $u,v\in V$. 
Note that $uv\in E$ if and only if $vu\in E$. 
Let $\T := \{z\in \mathbb{C} ~:~ |z| = 1\}$ be the multiplicative group of unimodular complex numbers. 
Then the mapping $\psi: E \mapsto \T$, with $\psi(uv)=\psi(vu)^{-1}$, is called a \textit{gain function}, and the tuple $\Psi = (G,\psi)$ is called a \textit{(unit) gain graph.}
The graphs $G$ are assumed to be connected throughout, though it should be noted that all gain graphs with exactly two distinct eigenvalues whose underlying graphs are not connected may be constructed by taking the disjoint union of two smaller such gain graphs whose distinct eigenvalues coincide. 

A cycle in a gain graph is said to be a (non-empty) circular walk in which the only repeated vertices are the first and the last.
For a given cycle, its gain is said to be the product of the gains of the traversed arcs. 
Mathematically, if the oriented cycle $C^\rightarrow$ is traversed by consecutively walking the arcs $u_1u_2,u_2u_3,...,u_{m-1}u_m,u_mu_1$, then the gain of $C^\rightarrow$ is given by $\phi(C^\rightarrow) = \psi(u_1u_2)\psi(u_2u_3)\cdots\psi(u_{m-1}u_m)\psi(u_mu_1)$. 
Note that if the cycle is traversed in reverse order, then $\phi(C^\leftarrow)=\overline{\phi(C^\rightarrow)}.$
Since the real part of the cycle gain typically contains all of the necessary information (see Theorem \ref{thm: coefs}), we will simply write $\phi(C)$, from now on.

For a given gain graph $\Psi$, its underlying graph is obtained by assigning every arc gain $1$. 
We define the \textit{underlying graph operator} $\Gamma(\cdot)$, that maps a gain graph to its underlying graph. 
A given graph $G$ is said to be \textit{$k$-regular} if every vertex has $k$ neighbors and \textit{bipartite} if it contains no odd-sized cycles.
Gain graphs are said to be $k$-regular and bipartite when their underlying graphs are. 

The main discussion in this work is concerned with the \textit{gain matrix} $A(\Psi)$ of $\Psi$, whose entries $[A(\Psi)]_{uv}$ are given by the corresponding $\psi(uv)$.
For a given gain graph $\Psi$ with gain matrix $A$, its \textit{characteristic polynomial} $\chi(\lambda)$, is said to be the characteristic polynomial of its gain matrix. That is, $\chi(\lambda) = \det\left(\lambda I-A\right).$ 
The \textit{eigenvalues} of $\Psi$ are the roots of $\chi(\lambda)$; the collection of eigenvalues is called the \textit{spectrum}, which is often denoted 
\[\Sigma_\Psi=\left\{\theta_1^{[m_1]},\theta_2^{[m_2]},\ldots,\theta_p^{[m_p]}\right\},\]
where $\theta_1,\ldots,\theta_p$ are the $p$ distinct eigenvalues of $\Psi$ and $m_1,\ldots,m_p$ are their respective multiplicities; not necessarily distinct eigenvalues are denoted by $\lambda_j$, $j\in[n]$. 
Note that the gain matrix is Hermitian, so it is diagonalizable with real eigenvalues. 

Two gain graphs $\Psi$ and $\Psi'$ are said to be \textit{isomorphic} (denoted $\Psi\cong\Psi'$) if they are equal, up to a relabeling of the vertices. 
Let $S$ be a diagonal matrix with diagonal elements in $\T$. Then $\Psi'$ is said to be obtained from $\Psi$ by a \textit{diagonal switching} if $A(\Psi') = S^{-1}A(\Psi)S.$
The \textit{converse} of $\Psi$ is obtained by inverting every edge gain of $\Psi$. 
\begin{definition}
\label{def: sw iso}
Two gain graphs are said to be \textit{switching isomorphic} if one may be obtained from the other by a sequence of diagonal switches, possibly followed by taking the converse and/or relabeling the vertices. 
Switching isomorphism of $\Psi$ and $\Psi'$ is denoted $\Psi\sim\Psi'$.
\end{definition}
If the vertices are not relabeled, the above is also known as \textit{switching equivalence}; we choose to slightly amend the definition to fully encompass the equivalence relation. 
A pair of convenient results when dealing with switching isomorphism are the following.
\begin{lemma}[\cite{reff2016oriented}, \cite{wissing2020spectral}]
\label{lemma: reff fundamental cycles}
Let $\Psi_1$ and $\Psi_2$ be gain graphs with the same underlying graph $G$. 
For every cycle $C$ in $G$ it holds that $\text{Re } \phi_1(C)=\text{Re }\phi_2(C)$ if and only if $\Psi_1$ and $\Psi_2$ are switching equivalent.  
\end{lemma}
As a minor note of care, we mention that the above does concerns switching \textit{equivalence}, rather than switching \textit{isomorphism.} 
Nevertheless, the above may be immediately applied to show that two gain graphs are not switching isomorphic if a one contains a cycle with gain $\phi$, while the other does not. (Note that diagonal switching does not affect cycle gains.)
\begin{lemma}[\cite{wissing2020spectral}]
\label{lemma: make a tree equal}
Let $G$ be a graph and let $\Psi_1=(G,\psi_1)$, $\Psi_2=(G,\psi_2)$. 
Let $T\subseteq E(G)$ be a spanning tree of $G$. 
Then there exists a diagonal switching matrix $X$ such that the gain graph $\Psi_2'=(G,\psi_2')$, obtained from $\Psi_2$ by $A(\Psi_2')=X^{-1}A(\Psi_2)X$, satisfies $\psi_1(e) = \psi'_2(e)$ for all $e\in T$. 
\end{lemma}
The above is most often applied to be able to assume the gains of a spanning tree in a gain graph, whose exact structure, subject to given constraints, is in question. 

We will sometimes be dealing with induced subgraphs. 
If $U\subseteq V$, then $G[U]$ denotes the graph that is obtained by removing all vertices in $V\setminus U$, and all edges (arcs) that are incident to said vertices.
Note that eigenvalue interlacing may then be applied, since the gain matrices are Hermitian.

A graph $G$ is called an elementary graph if each of its connected components is either an edge or a cycle.
The characteristic polynomial of a gain graph may be obtained from its elementary subgraphs as follows.

\begin{theorem}[\cite{samanta2019}]\label{thm: coefs}
Let $\Psi$ be a unit gain graph with underlying graph $G$. Let $\chi(\lambda)= \lambda^n + a_1\lambda^{n-1} + \ldots + a_n$ be the characteristic polynomial of $A(\Psi)$. Then
\[a_j = \sum_{H\in\mathcal{H}_i(G)}(-1)^{p(h)}2^{c(H)} \prod_{C\in\mathcal{C}(H)}\re{\phi(C)},\]
where $\mathcal{H}_i(G)$ is the set of all elementary subgraphs\footnote{Here, subgraphs may be vertex-induced, edge-induced or both.} of $G$ with $i$ vertices, $\mathcal{C}(H)$ denote the collection of all cycles in $H$, and $p(H)$ and $c(H)$ are the number of components and the number of cycles in $H$, respectively. 
\end{theorem}

Throughout, the identity matrix, the all-ones matrix and the zero matrix are denoted $I,J$ and $O$, respectively. 
Occasionally, a subscript is added to clarify its dimensions. 
An order-$n$ matrix is said to be \textit{cyclic} if its entries $c_{i,j}$ satisfy $c_{i+1,j+1}=c_{i,j}$ and $c_{i+1,1}=c_{i,n}$ for all $i,j\in[n-1]$, and $CM(\cdot)$ is the cyclic matrix operator that generates a square matrix based on its first row. 

Finally, we include a few often-used notions and their names.
As usual, a complete graph of order $n$ is denoted $K_n$, and a complete $k$-partite graph is denoted $K_{n_1,\ldots,n_k}$. 
Lastly, $C_n$ is the cycle of order $n$.
In the interest of brevity, we often denote by respectively $\varphi, \omega$ and $\gamma$ the third, sixth and eighth roots of unity. 
We conclude by offering an overview of the drawing conventions used in illustrations throughout, shown in Table \ref{tab: drawing conventions}. 
\begin{table}[h!]
    \centering
    \begin{tabular}{ll}
\toprule
Gain     & Drawing \\ \midrule
$\psi(uv)=1$ &       \begin{tikzpicture}
                        \node[vertex] (1) at (0,0) {$u$};
                        \node[vertex] (2) at (2,0) {$v$};
                        \draw[edge] (1) to node{} (2);
                        \end{tikzpicture}      \\
$\psi(uv)=-1$ &   \begin{tikzpicture}
                        \node[vertex] (1) at (0,0) {$u$};
                        \node[vertex] (2) at (2,0) {$v$};
                        \draw[negedge] (1) to node{} (2);
                        \end{tikzpicture}       \\
$\psi(uv)=\omega$  &    \begin{tikzpicture}
                    \node[vertex] (1) at (0,0) {$u$};
                    \node[vertex] (2) at (2,0) {$v$};
                    \draw[arc] (1) to node{} (2);
                    \end{tikzpicture}      \\
$\psi(uv)=-\omega = \varphi^2$  &    \begin{tikzpicture}
                        \node[vertex] (1) at (0,0) {$u$};
                        \node[vertex] (2) at (2,0) {$v$};
                        \draw[negarc] (1) to node{} (2);
                    \end{tikzpicture}      \\ 
$\psi(uv)=\gamma$  &    \begin{tikzpicture}
                        \node[vertex] (1) at (0,0) {$u$};
                        \node[vertex] (2) at (2,0) {$v$};
                        \draw[arc,dotted] (1) to node{} (2);
                    \end{tikzpicture}      \\
$\psi(uv)=1$, fixed ex ante  &    \begin{tikzpicture}
                        \node[vertex] (1) at (0,0) {$u$};
                        \node[vertex] (2) at (2,0) {$v$};
                        \draw[tree] (1) to node{} (2);
                    \end{tikzpicture}      \\ \bottomrule
\end{tabular}
    \caption{Drawing conventions}
    \label{tab: drawing conventions}
\end{table}

\subsection{Regularity in gain graphs}
Unlike for its graph analog, regularity of a gain graph $\Psi$ is not characterized by its spectrum.
To see this, \citet{belardo2019open} offer an example pair that are cospectral to one-another, while one is regular and the other is not. 
However, this pair does not feature two \textit{connected} signed graphs. 

It seems interesting to ask whether such examples may also be constructed under the assumption of connectedness. 
With relative ease, one finds a number of small examples to confirm the claim, even when both halves of the cospectral pair are required to be connected.
Figure \ref{fig: regularity not characterized connected version} illustrates one such example, though arbitrarily large ones may also be constructed. 
One such construction is provided below.
Let $K_{p,q}$ be the complete bipartite graph, whose nonzero eigenvalues are $\pm\sqrt{pq}$ with multiplicity $1$. 
Furthermore, let $K^*_{p,q,r}$ be a complete tripartite graph whose $3$-cycles all have gain $i$ and whose closed $4$-walks all have gain $1$; this graph has exactly two nonzero eigenvalues which are exactly $\pm \sqrt{pq+qr+rp}$.

\begin{figure}[t]
    \centering
    \begin{subfigure}[b]{.35\textwidth}\centering
    \begin{tikzpicture}
            \node[vertex] (1) at (1,.866) {};
            \node[vertex] (3) at (1,-.866) {};
            \node[vertex] (2) at (.5,0) {};
            \node[vertex] (5) at (-.5,0) {};
            \node[vertex] (6) at (-1,-.866) {};
            \node[vertex] (4) at (-1,.866) {};
            
            \draw[edge] (1) to node{} (2);
            \draw[edge] (1) to node{} (3);
            \draw[edge] (1) to node{} (4);
            \draw[edge] (2) to node{} (3);
            \draw[edge] (2) to node{} (5);
            \draw[edge] (3) to node{} (6);
            \draw[edge] (4) to node{} (5);
            \draw[edge] (4) to node{} (6);
            \draw[edge] (5) to node{} (6);
    \end{tikzpicture}
    \caption{}
    \end{subfigure}
    \begin{subfigure}[b]{.35\textwidth}\centering
    \begin{tikzpicture}
            \node[vertex] (1) at (-1,-.866) {};
            \node[vertex] (2) at (1,-.866) {};
            \node[vertex] (3) at (1,.866) {};
            \node[vertex] (4) at (-1,.866) {};
            \node[vertex] (5) at (-.5,0) {};
            \node[vertex] (6) at (.5,0) {};
            
            \draw[edge] (5) to node{} (1);
            \draw[edge] (5) to node{} (2);
            \draw[edge] (5) to node{} (3);
            \draw[edge] (5) to node{} (4);
            \draw[edge] (5) to node{} (6);
            \draw[edge] (6) to node{} (1);
            \draw[edge] (6) to node{} (4);
            \draw[edge] (6) to node{} (2);
            \draw[edge,dashed] (6) to node{} (3);
    \end{tikzpicture}
    \caption{}
    \end{subfigure}
    \caption{Cospectral connected signed graphs, of which only (a) is regular.}
    \label{fig: regularity not characterized connected version}
\end{figure}
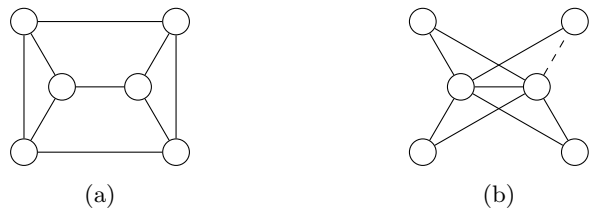
\begin{proposition}[\citet{xu2020complex}]\label{prop: rank 2}
Let $\Gamma$ be a connected graph and let $\Psi=(\Gamma,\psi)$ be a gain graph. 
Then Rank$(\Psi)=2$ if and only if either $\Psi\sim K_{p,q}$ or $\Psi\sim K_{p,q,r}^*.$
\end{proposition}

By using their respective spectra, as above, one finds the following.
\begin{corollary}\label{cor: cospec Kmm}
Let $m,j,s,t$ be natural numbers such that $m=j(s^2 + st + t^2)$. 
Then $K_{m,m}$ is cospectral to $K^*_{js^2, jt^2, j(s+t)^2}$.
\end{corollary}
The above shows a construction of arbitrarily large pairs of connected gain graphs, of which exactly one (i.e., $K_{m,m}$) is regular, while the other not necessarily is.
Note that this construction generalizes a remark that first appeared in \cite{mohar2016hermitian}.

As an aside, we note that one may characterize all gain graphs with rank $3$ with a straightforward, though tedious, forbidden subgraph approach. 
This collection may loosely be described as all gain graphs switching equivalent to a twin expansion \cite{wissing2019negative} of a triangle (not $K^*_{1,1,1}$) or a rank-$3$ gain graph on $K_4$. 
This has been proven by the authors; the details are omitted.

\subsection{Two distinct eigenvalues}
\label{subsec: two ev}
The main body of this work is concerned with gain graphs whose gain matrix has precisely two distinct eigenvalues; we will commonly call such objects \textit{two-eigenvalue gain graphs}.
Suppose that the two-eigenvalue gain graph $\Psi$ has eigenvalues $\theta_1$ and $\theta_2$ with multiplicities $m$ and $n-m$, respectively. 
If $a=\theta_1 + \theta_2$ and $k=-\theta_1\theta_2$, then the gain matrix $A$ of $\Psi$ satisfies
\begin{equation}\label{eq: minimal polynomial} A^2 = aA + kI.\end{equation}
This implies that $\Psi$ is $k$-regular, and hence that $k$ is integer. 
Moreover, since clearly $k>0$ (otherwise $\Psi$ is the empty graph), $A$ has full rank. 
Additionally, since $A$ and $-A$ have opposite eigenvalues, we will consider without loss of generality only the case with $a\geq 0$, which since $\tr A=0$ implies $m\leq n/2$.
Also note that $a\leq n-2$ 
with equality if and only if $\Psi$ is switching equivalent to a complete graph, with distinct eigenvalues $k$ and $-1$.

It is not hard to see that the eigenvalues of $\Psi$ are the square roots of rational numbers. 
Indeed, using that $0 = \tr A  = m\theta_1 + (n-m)\theta_2$ and $nk=\tr A^2 = m\theta_1^2 + (n-m)\theta_2^2$, it follows that 
\[\theta_1 = \sqrt{\frac{k(n-m)}{m}}~~\text{and}~~\theta_2 = -\sqrt{\frac{km}{n-m}}.\]
Moreover, by applying the quadratic formula to \eqref{eq: minimal polynomial}, we also have 
\[\theta_1 = \frac{a + \sqrt{a^2 + 4k}}{2}~~\text{and}~~\theta_2 = \frac{a - \sqrt{a^2 + 4k}}{2}.\]
If $a$ is integer, the following result from \citet{ramezani2018constructing} carries over. 
\begin{lemma}
\label{lemma: perfect square}
Let $\Psi$ be a two-eigenvalue gain graph, and let $a\in \mathbb{Z}$. Then either 
\begin{enumerate}
    \item [(i)] $a=0$ and the eigenvalues are $\pm \sqrt{k}$, or
    \item [(ii)] $a\not=0$ and $a^2 + 4k$ is a perfect square,
\end{enumerate}
\end{lemma}
\begin{proof}
The first part follows by plugging in $a=0$ into \eqref{eq: minimal polynomial}. 
The second part is shown by contradiction. 
Suppose that $a^2+4k$ is not a perfect square, so that $\theta_1$ is irrational. 
Then, since the characteristic polynomial of $\Psi$ is a monic integral polynomial, the algebraic conjugate of $\theta_1$, i.e., $\theta_2$, occurs as an eigenvalue of $A$ with the same multiplicity; say $m$. 
Now, since the trace of $A$ equals zero, it follows that 
$m\theta_1 + m\theta_2=ma = 0,$
and hence $a=0$, contradiction. 
\end{proof}
Equivalently, one may formulate this in terms of $n,m$ and $k$. 
\begin{lemma}\label{lemma: perfect square noninteger a}
Let $\Psi$ be a two-eigenvalue gain graph. 
If $a\in \mathbb{N}$ then $\frac{kn^2}{m(n-m)}$ is a perfect square.
\end{lemma}
\begin{proof}
If $a>0$, then by Lemma \ref{lemma: perfect square}, $a^2+4k=b^2$, for some $b\in\mathbb{N}$. 
We may rewrite to obtain:
\[~b^2=a^2+4k= \left( \sqrt{k(n-m)m^{-1}} -k\left(\sqrt{k(n-m)m^{-1}}\right)^{-1} \right)^2+4k= \frac{kn^2}{m(n-m)}.\]
\end{proof}

However, contrary to \cite{ramezani2018constructing}, the current context does not guarantee that $a$ is integer. 
Consider the following example, constructed from an equiangular tight frame of $7$ vectors in dimension $3$, which is closely related to the Fano plane. 
\example{edwin counterexample 1}{
Let $A=\frac14\sqrt{2}(I-J-i\sqrt{7}(N-N^{\top}))$, where $N=CM\left(\begin{bmatrix}0 &1& 1& 0& 1& 0& 0\end{bmatrix}\right).$
Then $a=\frac{1}{2}\sqrt{2}$ and the eigenvalues of $A$ are $2\sqrt{2}$ and $-\frac32\sqrt{2}$.
Moreover, note that $a^2 + 4k = 24\frac12,$ which is not a perfect square. 
}
The construction in Example \ref{edwin counterexample 1} is a member of an infinite family of gain graphs with two eigenvalues, that is based on a particular tight frame. 
Details regarding said family can be found in Section \ref{subsec: tight frames}.

As we will see over the course of this work, gain graphs that exhibit the desired extreme spectral behavior are extremely rare. 
This is of particular interest in the discussion concerning spectral characterizations of gain graphs. 
A gain graph is said to be (weakly) determined by its spectrum if any other gain graph with the same spectrum is switching isomorphic to the first. 
Indeed, since most of the here obtained examples are only examples of a given order with a particular spectrum, the above property oftentimes follows easily. 
The foremost obstacle to overcome is the possibility that disjoint unions of smaller graphs coincide spectrally. 
We will not point out every instance of spectral determination under connectedness, except in the summarizing Theorems \ref{thm: mult at most 3} and \ref{thm: summary}.

\section{Constructions}                
\label{sec: constructions}    
Somewhat unsurprisingly, there are various fairly well understood areas that are linked to the here considered notion. 
In this section, we will showcase these links and build on existing theory to obtain various two-eigenvalue gain graphs.

\subsection{Weighing matrices}
\label{subsec: weighing}

A complex unit weighing matrix is an $n\times n$ matrix $W$ with entries in $\T$ such that $WW^*=kI$, for some $k$. 
Real weighing matrices have been quite extensively studied (see \cite{harada2012classification}), and their complex generalizations have recently been getting more and more attention, too. 
For example, \citet{best2013unit} characterized all complex unit weighing matrices (simply \textit{weighing matrix,} hereafter) with weight at most $4.$ 

Note that since a weighing matrix is square and $WW^*=W^*W=kI$, a Hermitian weighing matrix $W$ with a zero diagonal characterizes a unit gain graph with eigenvalues $\pm\sqrt{k}$.  
The smallest nontrivial example of this is 
\[W_4=\begin{bmatrix}0&1 &1& 1 \\ 1 &0 &i &-i\\ 1& -i &0&i \\ 1& i &-i&0 \end{bmatrix}.\]
Following a convention for Hadamard matrices \cite{abiad2019graph}, weighing matrices are said to be \textit{graphical} when they are Hermitian and their diagonal is constant. 
This may only occur when the constant diagonal has value $\delta\in \{0,1,-1\}$; the corresponding gain graph is then obtained as $W-\delta I$.
By construction, such gain graphs have distinct eigenvalues $-\delta\pm\sqrt{k}.$ 

Below, we will mainly consider the generic case with $\delta=0$. 
It should be noted that similar considerations are possible when $\delta\not=0$, although these cases are considerably more restrictive. 
Indeed, if $\delta\not=0$ then the adjacency matrix $A=W-\delta I$ satisfies $A^2 = -2\delta A + (k-1)I$, where $W^2=kI$, and thus $4k$ must be a perfect square, by Lemma \ref{lemma: perfect square}. (Note that $A$ has degree $k-1$.)
An example of such a case is the complete graph $K_4$, whose adjacency matrix $A$ is related to the graphical Hadamard matrix $A-I$.

Another interesting link to the field of weighing matrices appears when one restricts oneself to the class of bipartite gain graphs. 
Indeed, if $\Psi$ is bipartite, then its gain matrix $A$ may be written as 
\begin{equation}A=\begin{bmatrix} O & B\\B^* & O\end{bmatrix},\label{eq: bipartite gain graph weighing}\end{equation}
and thus\footnote{Note that $a=0$ follows since bipartite gain graphs have symmetric spectra.}
\[A^2 = \begin{bmatrix} B^*B & O \\ O & BB^*\end{bmatrix}=kI \iff BB^*=B^*B=kI.\]
That is, $\Psi$ has exactly two distinct eigenvalues if and only if $B$ is a square
weighing matrix of weight $k$. 
In our below classification, we will will denote the bipartite gain graph obtained from a weighing matrix $B$ as in \eqref{eq: bipartite gain graph weighing} by $IG(B)$. 

We have to place a note of care here. 
It follows easily that the direct sum of any two weighing matrices of equal weight is again a weighing matrix. 
However, by the same token, if $W$ is the disjoint union of $W'$ and $W'',$ then $IG(W)$ is disconnected.
Since we assume connectedness throughout, we require that $W$ is irreducible.

The smallest nontrivial examples of unit weighing matrices are 
\[W_2=\begin{bmatrix}1 &1 \\ 1 &-1 \end{bmatrix}, ~W_3=\begin{bmatrix}1 &1& 1 \\ 1 &\f &\ff \\ 1& \ff &\f\end{bmatrix} \text{~and~} W_4.\]
In fact, we may draw some additional conclusions, regarding on these matrices. 
\begin{proposition}\label{prop: deg 2 characterization}
Let $\Psi$ be a connected, order-$n$ gain graph with eigenvalues $\pm\sqrt{2}$. 
Then $\Psi$ is switching isomorphic to $IG(W_2)$.
\end{proposition}
\begin{proof}
Immediate from \cite[Thm. 10]{best2013unit}, since any irreducible weighing matrix of weight $2$ is equivalent to $W_2$. 
\end{proof}
Using the subsequent result from \citet{best2013unit} that characterizes weighing matrices of weight $3$, one also readily finds the following, analogous result.  
\begin{proposition}\label{prop: deg 3 2evggs}
Let $\Psi$ be a connected, order-$n$ gain graph with eigenvalues $\pm \sqrt{3}$. 
Then $\Psi$ is switching isomorphic to either $W_4$, $IG(W_3)$ or $IG(W_4)$. 
\end{proposition}
In fact, it turns out that these graphs (and $K_4$) are the only two-eigenvalue gain graphs with degree $3$; this is shown formally in Section \ref{sec: small degree}. 
Further examples of  weighing matrices include \[W_5=CM\left(\begin{bmatrix}0& 1& \f& \f& 1\end{bmatrix}\right)\text{~and~}W_7=CM\left(\begin{bmatrix}-1 &1& 1& 0& 1& 0& 0\end{bmatrix}\right).\]
Of course, many more examples of weighing matrices may be (and have been) constructed, though we will not explicitly list them here. 
Several method to generate such examples are discussed in the next section. 

It should be noted that cubelike graphs, in particular, may often be equipped with a gain function such that the corresponding gain matrix is a weighing matrix. 
For example, $IG(W_4)$ is (switching equivalent to) a signed cube and by taking Kronecker products of the $2\times 2$ Pauli matrices, \citet{alon2020unitary} construct such gain graphs on the folded $k$-cube, and certain other Cayley graphs on $\mathbb{Z}_2^{k-1}$. 

\subsection{Recursive constructions}
\label{subsec: recursive doubles}
In the previous section, we have seen a construction that takes a given weighing matrix, and turns it into a gain graph with exactly two eigenvalues. 
In fact, if such a weighing matrix $W$ has a zero diagonal, such that it characterizes a gain graph $\Psi$, then $IG(W)$ in a sense doubles $\Psi.$
The following was effectively proven above, below \eqref{eq: bipartite gain graph weighing}.
\begin{lemma}
Let $\Psi$ be an order-$n$ gain graph with exactly two eigenvalues $\pm\sqrt{k}$. Then $IG(\Psi)$ has order $2n$ and eigenvalues $\pm\sqrt{k}$. 
\end{lemma}
One might wonder whether the reverse also holds, when the trivial counterexamples such as $W_2$ and $W_3$ are excluded. 
This is not true, as is shown in the example below.\footnote{Note that $B$ is equivalent to $W_4$ under the operations listed in \cite{best2013unit}. However, $B$ is not graphical, while $W_4$ is.} 
\example{ex: IG counterexample}{
Let $B$ be the matrix defined as
\[B=\begin{bmatrix}1 & 1 & 1 & 0\\ 1 & -1 & 0 & 1\\ 1 & 0 & -1 & -1\\ 0 & 1 & -1 & 1\end{bmatrix}.\]
Then the signed cube $IG(B)$ has eigenvalues $\pm\sqrt{3}$, while $B$ itself is not graphical. 
}
Many such examples may be constructed, by using a set of operations that map a given weighing matrix to another, which does not necessarily preserve the Hermitian property, or the zero diagonal. 
These operations include permuting the rows (eq. columns) or multiplying a row (eq. column) by a number in $\T$; see \cite{best2013unit} for details. 

An idea similar to the doubling operation above was recently used by \citet{huang2019induced} in his proof of the Sensitivity Conjecture of Nisan and Szegedy on Boolean functions. 
For a given Hermitian matrix $W$ with exactly two distinct eigenvalues $\pm\sqrt{k}$, one easily finds that 
\begin{equation}
A=\begin{bmatrix}W & I\\ I & -W\end{bmatrix}
\end{equation}
has distinct eigenvalues $\pm\sqrt{k+1}$. 
In particular, this construction was used by Huang to construct signed $n$-cubes; see \cite{belardo2019open} for more info. 
We will call this construction Huang's Negative Double, denoted $ND(\Psi)$. 

Stani\'{c} \cite{stanic2020decomposition} observed that under the same conditions,  Sylvester's recursive construction for Hadamard matrices carries over to the current paradigm. That is, the matrix
\begin{equation} B=\begin{bmatrix} W &W \\ W &-W\end{bmatrix}\end{equation} 
also has two distinct eigenvalues, i.e., $\pm \sqrt{2k}$.
This construction is, in turn, called the Sylvester Double and denoted $SD(\Psi)$.
Moreover, we obtain variation on the above by adding an identity component to the off-diagonal blocks. 
Specifically, the matrix 
\begin{equation} B=\begin{bmatrix} W &W+iI \\ W-iI &-W\end{bmatrix}\label{eq: stanic variation}
\end{equation}
has eigenvalues $\pm \sqrt{2k+1}$.
This operation will be denoted $SD^*(\Psi)$, hereafter. 

It should be noted that the above results hold for general weighing matrices $W$, even though we are technically only interested in those with unit entries and a zero diagonal. 
Observe that if $W$ admits to these conditions, then so does each of its doubles, described above; thus, said constructions are directly applicable to gain graphs. 
Formally, we have the following. 
\begin{lemma}
Let $\Psi$ be an order-$n$ gain graph with distinct eigenvalues $\pm \sqrt{k}$. 
Then the distinct eigenvalues of $ND(\Psi)$ are $\pm\sqrt{k+1}$ and the distinct eigenvalues of $SD(\Psi)$ are $\pm\sqrt{2k}.$
Finally, the distinct eigenvalues of $SD^*(\Psi)$ are $\pm\sqrt{2k+1}$.  
\end{lemma}

A final construction that follows a similar pattern was provided by \citet{greaves2012cyclotomic} in his classification of cyclotomic matrices over the Gaussian and Eisenstein integers. 
A concrete description of these is as as follows. Let $C$ be an order-$t$ weighing matrix with weight 1 and a zero diagonal, such that $C+C^*$ is also a unimodular matrix. 
Then, construct $A$ as
\begin{equation}
    A=\begin{bmatrix}C+C^* & C-C^*\\ C^*-C  &-C-C^{*}\end{bmatrix} \label{eq: greaves construction}.
\end{equation}
One is easily convinced that $A^2=4I$, and thus $A$ has eigenvalues $\pm 2$. 

Note that $C+C^*$ is the gain matrix of a cycle with gain $x$, for some $x\in \T$. 
Using Lemma \ref{lemma: make a tree equal}, we may without loss of generality switch such that all but one entry of $C$ equal one, such that the final entry equals $x$. 
For $C$ defined in such a way, and $A$ obtained from $C$ as in \eqref{eq: greaves construction}, we say that $A$ is the gain matrix of a \textit{toral tesselation}  graph \cite{greaves2012cyclotomic}, which is denoted $T_{2t}^{(x)}$. 
Note that the graphs $T_{2t}^{(x)}$, $T_{2t}^{(\overline{x})}$, $T_{2t}^{(-x)}$, and $T_{2t}^{(-\overline{x})}$ are all switching isomorphic.

Finally, we note that a variation on \eqref{eq: greaves construction} similar to \eqref{eq: stanic variation} is also possible. 
That is, 
\[B=\begin{bmatrix}C+C^* & C-C^*+I\\ C^*-C+I  &-C-C^{*}\end{bmatrix},\]
which is later said to be a donut graph, has eigenvalues $\pm\sqrt{5}$. 
We discuss this construction in more detail in Section \ref{sec: degree 5}, when we use it to construct infinite families of gain graphs with eigenvalues $\pm\sqrt{5}$, for every even $n\geq 8$. 

To end this section, the authors would like to express some interest in similar recursive constructions that do not require its blocks to be weighing matrices. 
In particular, it seems plausible that gain graphs with two distinct eigenvalues that do not sum to zero, may also be expanded into larger graphs that keep much of their structure, and thereby have exactly two distinct eigenvalues, as well. 
However, such constructions are unknown to the authors, at the time of writing.

\section{Lines separated by few distinct angles}                  
\label{sec: systems of lines}
Interestingly, the matter at hand has various links to other well-studied fields that are more geometric and algebraic in nature.
In particular, numerous topics that are all based on of a system of lines in complex space are naturally tied to the highly structured matrices that we are interested in. 
The connection between systems of lines and highly symmetric graphs has been explored before, most notably in the classification of graphs with least eigenvalue $-2$ by \citet{cameron1991line}, which was recently extended to signed graphs by \citet{greaves2015edge}. 
In this work, we explore a similar connection in a much more general setting. 

In this section, we will touch on several interconnected research areas that are concerned with these peculiar systems of lines. 
These topics include (Tight) Frames \cite{waldron2018introduction}, Mutually Unbiased Bases (MUB) \cite{durt2010mutually} and Symmetric Informationally Complete Positive Operator-Valued Measurements (SIC-POVM) \cite{geng2020minimal}. 
As we will see shortly, the more algebraic approach feeds naturally into a perspective that bounds the multiplicities, which leads to some beautiful examples of two-eigenvalue gain graphs, that are derived from well-known mathematical objects, such as the Witting polytope and the Coxeter-Todd Lattice.

\subsection{The Eisenstein matrix}
In an earlier work \cite{wissing2020spectral}, the authors have considered a Hermitian adjacency matrix for Signed Directed Graphs.
This matrix, which was called the Eisenstein matrix, after the group of unit Eisenstein integers $\T_6$ that make up its nonzero entries, may simply be considered to be the gain matrix of a gain graph. 
The current line of questioning does therefore apply. 
In this section, we will offer a brief intermezzo in which we will restrict the allowed edge gains to the entries of $\T_6$, in order to illustrate the perspective one might obtain by considering systems of lines. 

The attentive reader may have observed that almost all examples (excluding Example 1) have had either $a=0$ or $a=k-1$.
(Recall that $a=\theta_1+\theta_2$ and $k=-\theta_1\theta_2$, where $k$ is the degree of the corresponding gain graph.)
As an illustrative exercise, let us attempt to construct signed digraphs with exactly two distinct eigenvalues, such that $0<a<k-1$. 

An important detail to note here, is that $a\in\mathbb{Z}$ when the edge gains are restricted to $\T_6$. 
Indeed, since $A^2 = aA + kI$,
we have
\begin{equation}a = \sum_h A_{ih}A_{hj}A_{ji},\label{eq: sum of elements t6}\end{equation}
for some nonzero $A_{ij}$. 
In this particular case, \eqref{eq: sum of elements t6} then means $a$ is a sum of elements in $\T_6$. 
Moreover, since $a=\theta_1+\theta_2$, $a$ is real and it follows naturally that $a\in\mathbb{Z}$.
Hence, Lemma \ref{lemma: perfect square} may be applied.\footnote{The same conclusion can be reached by using that $A$ has a characteristic polynomial with integer coefficients, see e.g. \cite[Thm. 3.3]{wissing2020spectral}. It then follows that its minimal polynomial $\lambda^2-a\lambda-k$ has integer coefficients, as well.}
(As an aside, a parallel argument holds when the gains are restricted to $\T_4$.)

Since Lemma \ref{lemma: perfect square} applies, the tuple $(a,k)$ must satisfy $a^2+4k=b$ for some integer $b$.
One is easily convinced that the smallest value of $k$ such that $0<a<k-1$ and the above holds is $k=6$, in which case $a=1$. 
It then follows that $\theta_1=3$ and $\theta_2=-2$, and thus $m=2n/5$ and $n-m=3n/5$. 
This, in turn, implies that $n$ must be a multiple of $5$. 
In the below, we will consider a number of possible values $n$, and discuss possible examples of signed digraphs with two distinct eigenvalues and the before mentioned parameters. 

\begin{figure}[h!]
    \centering
    \begin{subfigure}[b]{.49\textwidth}\centering
    \begin{tikzpicture}
     \foreach \a in {1,2,...,5}{
     \node[vertex] (\a) at (18+\a*360/5: 2 cm) {};
     }
     \foreach \a in {6,7,...,10}{
     \node[vertex] (\a) at (18+36+\a*360/5: 2 cm) {};
     }     
     \foreach \a in {1,2,...,5}{
     \node (e\a) at (18+\a*360/5: 2.8 cm) {};
     }
     \foreach \a in {6,7,...,10}{
     \node (e\a) at (18+36+\a*360/5: 2.8 cm) {};
     }

\def\ebr{26}
\draw[negedge, bend left=15] (1) to node{} (2);
\draw[negedge, bend left=15] (2) to node{} (3);
\draw[negedge, bend left=15] (3) to node{} (4);
\draw[negedge, bend left=15] (4) to node{} (5);
\draw[negedge, bend left=15] (5) to node{} (1);

\draw[negedge, bend right=10] (6) to node{} (8);
\draw[negedge, bend right=10] (8) to node{} (10);
\draw[negedge, bend right=10] (10) to node{} (7);
\draw[negedge, bend right=10] (7) to node{} (9);
\draw[negedge, bend right=10] (9) to node{} (6);

\draw[negedge, bend right=5] (1) to node{} (6);
\draw[negedge, bend right=5] (6) to node{} (2);
\draw[negedge, bend right=5] (2) to node{} (7);
\draw[negedge, bend right=5] (7) to node{} (3);
\draw[negedge, bend right=5] (3) to node{} (8);
\draw[negedge, bend right=5] (8) to node{} (4);
\draw[negedge, bend right=5] (4) to node{} (9);
\draw[negedge, bend right=5] (9) to node{} (5);
\draw[negedge, bend right=5] (5) to node{} (10);
\draw[negedge, bend right=5] (10) to node{} (1);

\draw[edge, bend right=\ebr] (1) to (e6) to (e2) to (7);
\draw[edge, bend right=\ebr] (6) to (e2) to (e7) to (3);
\draw[edge, bend right=\ebr] (2) to (e7) to (e3) to (8);
\draw[edge, bend right=\ebr] (7) to (e3) to (e8) to (4);
\draw[edge, bend right=\ebr] (3) to (e8) to (e4) to (9);
\draw[edge, bend right=\ebr] (8) to (e4) to (e9) to (5);
\draw[edge, bend right=\ebr] (4) to (e9) to (e5) to (10);
\draw[edge, bend right=\ebr] (9) to (e5) to (e10) to (1);
\draw[edge, bend right=\ebr] (5) to (e10) to (e1) to (6);
\draw[edge, bend right=\ebr] (10) to (e1) to (e6) to (2);

    \end{tikzpicture}
    \caption{}\label{fig: ramezani example}
    \end{subfigure}
    \begin{subfigure}[b]{.49\textwidth}\centering
    \begin{tikzpicture}[scale=1]
    \def\offset{0}
     \foreach \a in {1,2,...,5}{
     \node[vertex] (\a) at (90-72+\a*360/5: 3 cm) {};
     }
     \foreach \a in {6,7,...,10}{
     \node[vertex] (\a) at (\offset+270+72+72+\a*360/5: 2.4 cm) {};
     }
     \foreach \a in {11,12,...,15}{
     \node[vertex] (\a) at (-\offset+270-72+\a*360/5: 1.2 cm) {};
     }

\draw[edge, bend right=25] (1) to node{} (2);
\draw[edge, bend left=25] (1) to node{} (5);
\draw[negarc, bend right=25] (2) to node{} (3);
\draw[negarc, bend right=25] (3) to node{} (4);
\draw[negarc, bend right=25] (4) to node{} (5);

\draw[negarc] (6) to node{} (1);
\draw[edge] (1) to node{} (8);
\draw[edge] (1) to node{} (10);
\draw[negarc] (6) to node{} (2);
\draw[negarc] (2) to node{} (7);
\draw[edge] (2) to node{} (9);
\draw[negarc] (7) to node{} (3);
\draw[negarc] (8) to node{} (3);
\draw[edge] (3) to node{} (10);
\draw[edge] (4) to node{} (6);
\draw[edge] (4) to node{} (8);
\draw[edge] (4) to node{} (9);
\draw[edge] (5) to node{} (7);
\draw[negarc] (5) to node{} (9);
\draw[edge] (5) to node{} (10);

\draw[edge] (1) to node{} (11);
\draw[edge] (2) to node{} (12);
\draw[edge] (3) to node{} (13);
\draw[edge] (4) to node{} (14);
\draw[edge] (5) to node{} (15);

\draw[negarc, bend left=25] (6) to node{} (13);
\draw[edge, bend right=0] (6) to node{} (14);
\draw[edge, bend right=25] (6) to node{} (15);

\draw[negarc, bend left=25] (11) to node{} (7);
\draw[negarc, bend left=25] (7) to node{} (14);
\draw[edge, bend right=0] (7) to node{} (15);

\draw[edge, bend right=0] (8) to node{} (11);
\draw[negarc, bend left=25] (12) to node{} (8);
\draw[negarc, bend right=25] (15) to node{} (8);

\draw[negarc, bend right=25] (11) to node{} (9);
\draw[edge, bend right=0] (9) to node{} (12);
\draw[edge, bend right=25] (9) to node{} (13);

\draw[negarc, bend left=25] (10) to node{} (12);
\draw[edge, bend right=0] (10) to node{} (13);
\draw[edge, bend right=25] (10) to node{} (14);

\draw[negarc, bend left=15] (13) to node{} (11);
\draw[negarc, bend left=15] (11) to node{} (14);
\draw[negarc, bend right=15] (12) to node{} (14);
\draw[edge, bend left=15] (12) to node{} (15);
\draw[negarc, bend right=15] (13) to node{} (15);

     \foreach \a in {16,17,...,20}{
     \node[vertex,fill=white] (\a) at (-\offset+270-72+\a*360/5: 1.2 cm) {};
     }

    \end{tikzpicture}
    \caption{}\label{fig: gq}
    \end{subfigure}
    \caption{Signed digraphs with spectra $\left\{3^{[2n/5]},-2^{[3n/5]}\right\}$. (In (b), all three vertices hit by each straight line through the center of the picture are pairwise adjacent; all such edges have gain 1.)}
\end{figure}
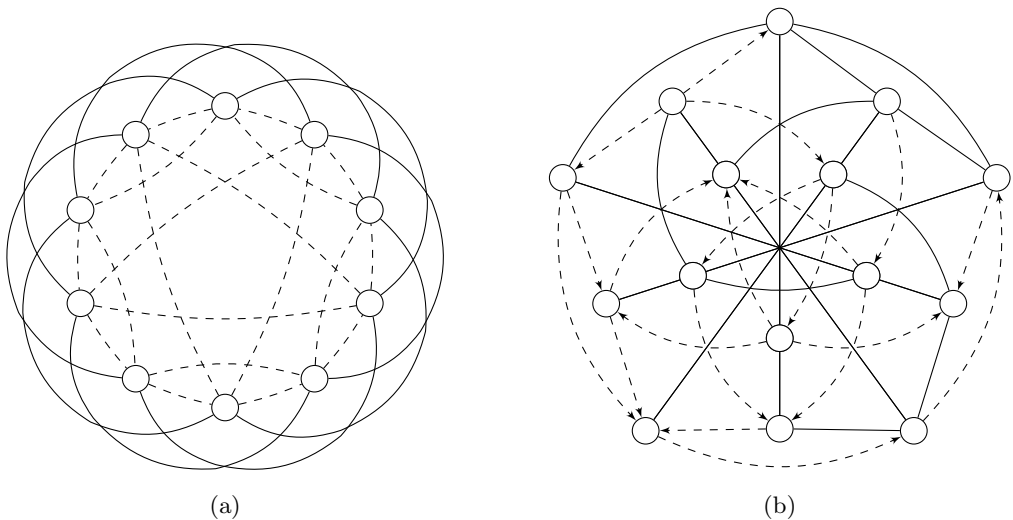

The smallest possible $n$ is $n=10.$ 
\citet{ramezani2018constructing} has constructed a signed graph (gains in $\T_2\subset \T_6$) with the above spectrum on the complement of the Petersen graph. 
Ramezani moreover shows that this example, illustrated in Figure \ref{fig: ramezani example}, is actually a member of an infinite family of signed graphs with two distinct eigenvalues on the triangular graphs\footnote{The triangular graph $\Delta(m)$ is the line graph of complete graph $K_m$.} $\Delta(m)$.

The next case, $n=15$, also admits an example with the above parameters that is related to the triangular graphs.
Below, we propose a construction on the generalized quadrangle $GQ(2,2)$, which is the complement of $\Delta(6)$, and illustrated in Figure \ref{fig: gq}.
To find the desired example, we employ the so-called hexacode \cite{conway2013sphere}: a 3-dimensional linear code of length $6$ over $GF(4)=\{0,1,\f,\ff\}$, where $\f$ and $\ff$ denote the third roots of unity. 
Specifically, the hexacode is defined by
\[H = \left\{\begin{bmatrix}p_2&p_1&p_0&f(1)&f(\f)&f(\ff)\end{bmatrix} ~:~ f(x)=p_2x^2+p_1x+p_0,~ p_2,p_1,p_0\in GF(4)\right\}.\]

In particular, the hexacode has $45$ elements, called codewords, of weight $4$, occurring in $15$ $1$-dimensional subspaces; i.e., lines through the origin. 
From each such subspace, we choose one nonzero codeword which we consider as a vector in $\mathbb{C}^6$. 
These vectors $v_h, h\in[15]$, represent our vertex set.
It is easily verified that each of the possible 15 supports occurs exactly once, and distinct supports can intersect in 2 or 3 positions. 
Moreover, the construction is such that if two codewords $v_h$ and $v_j$ have supports that intersect in 3 positions, then the corresponding inner product is always $v_h^*v_j=1+\f+\ff=0$. 
If two codewords have supports that intersect in 2 positions, then the inner product is either $1+1$, $\f+\f$, or $\ff+\ff$. Hence, we may define an Eisenstein matrix $\mathcal{E}$ (gain matrix) by
\[\mathcal{E}_{hj}=\frac12 v_h^*v_j,~h\not=j; \mathcal{E}_{jj}=0.\]
Note that the above indeed defines a signed directed graph, since $\mathcal{E}_{hj}\in \{0,1,\f,\ff\}$ and $\mathcal{E}$ is Hermitian. 
Moreover, if we set $M$ to be the matrix whose columns are the $v_j$, then according to the definition above $\mathcal{E}=\frac12 M^*M -2I,$ and thus
\begin{equation}\mathcal{E}^2 = \frac14 M^*MM^*M-2M^*M+4I = \frac{1}{2}M^*M +4I = \mathcal{E}+6I,\label{eq: gq22}\end{equation}
where the second equality follows since $MM^*=10I.$ 
(A similar fact should hold in general; this is formalized in Proposition \ref{prop: NN*=kI}, below.)
Lastly, note that by \eqref{eq: gq22}, it follows that $\mathcal{E}$ indeed has the desired spectrum $\{3^{[6]},-2^{[9]}\}.$ 

Note that, by definition, taking a different representative vector of a subspace will lead to a signed digraph that is switching equivalent with the original one. 
The corresponding equivalence class is, in fact, the only one with the desired spectrum with this particular underlying graph, thus yielding an easy spectral characterization. 
The particulars to this fact are quite tedious, and have been verified by computer. 

The construction above has some ties to previously studied objects. 
It is, for example, closely related to the so-called tilde-geometry \cite{pasini2002some}. 
Moreover, Figure \ref{fig: gq} is in a sense a quotient of the distance-regular antipodal $3$-cover of the collinearity graph of the generalized quadrangle of order $2$ \cite[p.~398]{brouwer1989distance}.
This, in turn, is a distance-regular graph that is defined on the above mentioned $45$ codewords of weight $4$, with adjacency of vertices $h$ and $j$ if 
$\mathcal{E}_{hj}=\frac12 v_h^*v_j$ equals $1$.

Most importantly, the above discussion sheds some  light on the way in which systems of lines ($1$-dimensional subspaces) are connected to the remarkable gain graphs that are the topic of this work.
In essence, one needs systems of lines that are either orthogonal to one-another, or are all separated by the same specific angle. Correspondingly, Ramezani's examples on $\Delta(m)$ can be described by the lines through the  vectors $e_h - e_j$ ($h<j$) in $\mathbb{R}^m$, where $e_h$ is a standard basis vector.
Systems of lines that are pairwise separated by the same angle (so-called \textit{equiangular} lines) will be of special interest, later in this section. 
First, we offer a little more general insight based on the above.

\subsection{Decomposition as a $\{0,\alpha\}$-set}
In the above, we have showcased a clear parallel between gain graphs with few distinct eigenvalues, and systems of lines that are separated by few distinct angles. 
The equation $\mathcal{E}=\frac{1}{2}M^*M-2I$ is particularly reminiscent of a Gram matrix, though it does need a little additional work. 

In general, if $A:=A(\Psi)$ for some gain graph $\Psi$ with smallest eigenvalue $\theta_{\min}$, whose multiplicity is $n-m$, then $I-\theta_{\min}^{-1}A$ is a positive semi-definite matrix with rank $m$. 
It can therefore be represented as the Gram matrix of (Hermitian) inner products of a set of $n$ unit vectors $\{u_1,u_2,\dots,u_n\}$ in complex space $\mathbb{C}^m$.
As before, the absolute values $|u_i^*u_j|$ of these inner products represent the angles between the lines through the unit vectors.
In our case, there are exactly 1 or 2 such angles.
Correspondingly, the inner product of every two distinct unit vectors has absolute value either zero or $-\theta_{\min}^{-1}$.
In the study of lines in (complex) space, this phenomenon is known as a $\{0,\alpha\}$-set, where $\alpha=-\theta_{\min}^{-1}$ is the non-orthogonal separation angle. 
\example{ex: gram}{
Recall $W_2 = \begin{bmatrix}1 & 1 \\1 & -1\end{bmatrix}$ and let $\Psi=IG(W_2)$. 
Then $\Psi$ has eigenvalues $\pm\sqrt{2},$ and thus
\[I-\theta_2^{-1}A(\Psi) = \begin{bmatrix} 1 & 0 & \frac12\sqrt{2}& \frac12\sqrt{2}\\
                                          0 & 1 & \frac12\sqrt{2} & -\frac12\sqrt{2}\\
                                          \frac12\sqrt{2} & \frac12\sqrt{2} & 1 & 0\\
                                          \frac12\sqrt{2} & -\frac12\sqrt{2} & 0 & 1\end{bmatrix},\]
which is the Gram matrix of the unit vectors
\[u_1=\begin{bmatrix} 1 & 0 \end{bmatrix}^\top,~~ 
  u_2=\begin{bmatrix} 0 & 1 \end{bmatrix}^\top,~~ 
  u_3=\frac12\sqrt{2}\begin{bmatrix} 1 & 1 \end{bmatrix}^\top,~~ 
  u_4=\frac12\sqrt{2}\begin{bmatrix} 1 & -1 \end{bmatrix}^\top.  \]
}
As before, the vector notation can be translated to a matrix (outer)product.
Let $N$ be the matrix whose columns are the vectors $u_1,\ldots,u_n$. 
Then, by the above, $N^*N = I-\theta_{\min}^{-1}A.$ 
However, one cannot carelessly expect that any $\{0,\alpha\}$-set of unit vectors will yield a two-eigenvalue gain graph, as is showcased in the following example. 

\example{ex: reverse}{
Let $\{u_1,\ldots,u_4\}$ be the collection of vectors
\[\left\{ \begin{bmatrix} 1 & 0 & 0 \end{bmatrix}^\top, 
\begin{bmatrix} \frac12 & \frac12\sqrt{3} & 0 \end{bmatrix}^\top, 
\begin{bmatrix} 0 & \frac13\sqrt{3}& \frac13\sqrt{6} \end{bmatrix}^\top, 
\begin{bmatrix} \frac12 & -\frac16\sqrt{3} & \frac13\sqrt{6} \end{bmatrix}^\top \right\}\]
and let $N = \begin{bmatrix} u_1& u_2& u_3& u_4\end{bmatrix}$. 
Then 
\[N^*N = \frac{1}{2}\begin{bmatrix} 2 & 1 & 0 & 1 \\ 1 & 2 & 1 & 0 \\ 0 & 1 & 2 & 1\\ 1 & 0 & 1 & 2\end{bmatrix} = I-(-2)^{-1}A(C_4),\]
where $C_4$ is the undirected four-cycle, which has spectrum $\{-2,0^{[2]},2\}.$
}
This goes to show that $n$ unit vectors that are separated by one of two angles do not, in general, suffice to find a two-eigenvalue gain graph. 
However, the following interesting fact gives us an easy characterization. 
\begin{proposition}\label{prop: NN*=kI}
Let $A$ be an order-$n$ Hermitian matrix with least eigenvalue $\theta_\min\not=0$, whose multiplicity is $n-m$, $m>0$, and let $N\in\mathbb{C}^{m\times n}$ be such that $N^*N=I_n-\theta_\min^{-1}A.$
Then $A$ has exactly two distinct eigenvalues if and only if $NN^*=zI_m$ for some $z\in \mathbb{R}$. 
\end{proposition}
\begin{proof}
 $A$ has exactly two distinct eigenvalues if and only if $N^*N\in\mathbb{C}^{n\times n}$ has a spectrum given by
 $\left\{0^{[n-m]},1-\theta_2^{-1}\theta_1^{[m]}\right\}$.  
 Since $N^*N$ and $NN^*$ coincide on the nonzero eigenvalues, this is equivalent with $NN^*\in\mathbb{C}^{m\times m}$ having a single eigenvalue $1-\theta_1\theta_2^{-1}$ with multiplicity $m$.
 This, in turn, occurs if and only if $NN^*=(1-\theta_1\theta_2^{-1})I$. 
\end{proof}
It should be noted that if $A$ has a zero diagonal and constant-norm nonzero entries, then $z=1-\theta_1\theta_2^{-1}=n/m \in \mathbb{Q}$.  
Furthermore, if the columns of $N$ form a $\{0,\alpha\}$-set of unit vectors such that $NN^*=zI$, then either $N^*N=I$ and the columns of $N$ form an orthonormal basis ($z=1$, $n=m$, and the matrix $A=0$), or $NN^*=I+\alpha A$ and $A$ is a two-eigenvalue gain graph.

To conclude this section, we touch on some interesting facts. 
In case two of the unit vectors from $u_1,u_2,\dots,u_n$ are scalar multiples of each other, then their inner product is a unit. 
This implies that $\theta_{\min}=-1$ and therefore that the corresponding unit gain graph is switching isomorphic to a complete graph. 
It follows that $m=1$ and the $u_j$ are simply unit complex numbers.
This somewhat trivial case will be excluded in the classifications discussed in Section \ref{sec: small degree}. 

Moreover, one should also exercise some care when taking induced subgraphs. 
Instinctively one might be keen to claim that a taking a subset of a given system of lines that constitutes a two-eigenvalue gain graph will yield another one. 
However, similarly to was displayed in Example \ref{ex: reverse}, above, one cannot carelessly remove columns from $N$ without affecting the entries of $NN^*$, the latter of which must be a multiple of the identity. 
In Section \ref{subsec: mubs}, we will discuss a way to take two-eigenvalue subgraphs. 
\subsection{Bounds}
\label{subsec: bounds}
An intuitive question related to the matter of lines in complex space has to do with existence.
Specifically: how many lines can there be in an $m$ dimensional space, such that the angle between each pair is one of a given number of possible angles. 
In general, the classic result by \citet{delsarte1975bounds}, that bounds the number $n$ of distinct lines in $\mathbb{C}^m,$ whose separation angles are all contained in some collection $\mathcal{A}_s=\{\alpha_1,\alpha_2,\ldots,\alpha_s\}$ tells us that 
\begin{equation}\label{eq: absolute bound}
    n\leq 
    \binom{m+s-1}{m-1}\binom{m+s-1-\varepsilon}{m-1},
\end{equation}
where $\varepsilon=1$ if $0\in \mathcal{A}_s$ and zero otherwise. 
This has been called the absolute bound for systems of lines; for a particularly concise proof, the interested reader is referred to \cite{koornwinder1976note}. 
As should be evident from the previous section, we will exclusively concern ourselves with the cases $\mathcal{A}_1=\{\alpha\}$ and $\mathcal{A}_2=\{0,\alpha\}$, for $\alpha\in(0,1).$
In these cases, \eqref{eq: absolute bound} reduces to $n\leq m^2$ and $n\leq \frac12m^2(m+1),$ respectively. 
However, the latter bound may be sharpened under particular circumstances. 
\begin{proposition}
\label{prop: sharper bound A2 set}
Let $\Psi$ be an order-$n$ gain graph with spectrum $\{\theta_1^{[m]},\theta_2^{[n-m]}\}$, such that $\Gamma(\Psi)$ has degree $k$ and an eigenvalue $-\frac{km}{n-m}$ with multiplicity $m'\geq 0$. Then $n\leq m^2+m'$. 
\end{proposition}
\begin{proof}
Let $\{v_j\}_{j=1}^n$ be the corresponding system of vectors in $\mathbb{C}^m$, and consider the matrix $M$ 
defined by $M_{hj} = \tr v_hv^*_h
v_jv^*_j = |v^*_hv_j|^2$. Then $M = I + \frac{n-m}{km} B$, where $B$ is the adjacency matrix of $\Gamma(\Psi)$.
Consider now the linear transformation $T:\mathbb{C}^n \mapsto \mathbb{C}^{m \times m}$ defined by $T(x)=\sum_j x_jv_jv_j^*$. 
Then, by the rank-nullity theorem, one has $n=\dim \text{Range}(T) + \dim \ker(T) \leq m^2+m'$.
Here, it is used that $\dim\ker(T)\leq m'$ since $T(x)=0$ implies $Mx=0$ and thus $x$ is in the eigenspace of $B$ for eigenvalue $-\frac{km}{n-m}$. 
\end{proof}

 The above directly ties into some of the well-studied geometric objects that will be discussed shortly. 
For example, if $\Gamma=K_{t\times m},$ then $m'=n/m-1$ and the above reduces to $n\leq m(m+1)$.
These are precisely the underlying graph and the bound that occur in the case of Mutually Unbiased Bases, which is treated in Section \ref{subsec: mubs}. 
Furthermore, it should be clear from the proof above that in the case of equality, the projectors $v_jv_j^*$ span $\mathbb{C}^{m\times m}$. 
If, additionally, $m'=0$ then they form a basis of $\mathbb{C}^{m\times m}$. 
This corresponds to the absolute bound in the case of an $\mathcal{A}_1$-set, which is attained by a 'symmetric, informationally complete positive operator-valued measurement,' which are treated in Section \ref{subsec: sic povms}.

To conclude this section, we touch on the relative bound, which originally was a bound on the eigenvalues of Seidel matrices by \citet[Lemma 6.1]{vanlint1966equilateral}.
It bounds the angle $\alpha$ of a system of $n$ lines (represented by unit vectors $u_j$) in $\mathbb{C}^m$ with mutual angles at most $\alpha$.
Although this bound does not have any implications for our work, it is good to notice that equality in this bound leads to a two-eigenvalue graph. 
This is immediately clear from the proof\footnote{Indeed, note  that in the case of equality $\tr YY^*=0$ and thus $Y=\sum_ju_ju_j^*-\frac{n}{m}I=0$, which implies that $\sum_j u_ju_j^*$, which in our notation is equal to $NN^*$,  is a multiple of $I$, and hence Proposition \ref{prop: NN*=kI} applies. 
Therefore, we indeed have equality if and only if the corresponding gain graph has two distinct eigenvalues.}  by \citet[Prop.~10.6.3]{brouwer2011spectra}.

\subsection{Tight frames}
\label{subsec: tight frames}
Another school of thought is concerned with the notion of \textit{tight frames}. 
Boiled down to its essence, a tight frame is an over-complete collection of vectors that span some vector space; that is, the frame contains some deliberately redundant members. 
Under the right circumstances, this redundancy is actually an advantage. 
Since their conception, tight frames have found the most use in the development of \textit{wavelets,} since decomposition into tight frames rather allowed for far simpler representation than orthonormal bases would. 

Let us provide the definitions \cite{waldron2018introduction}.
Formally, a frame is a set of vectors $\{v_k\}_{k\in\mathcal{K}}$ in a Hilbert space $\mathcal{H}$, indexed by some collection $\mathcal{K}$, that satisfy
\[c_1\|u\|^2\leq \sum_{k\in\mathcal{K}}|\langle u,v_k\rangle|^2\leq c_2\|u\| ~~\forall u\in \mathcal{H},\]
for some constants $0< c_1\leq c_2<\infty.$
In case $c_1=c_2$, then $\{v_k\}_{k\in\mathcal{K}}$ is said to be a \textit{tight frame.}
Moreover, in case $c_1=c_2=1$, then $u=\sum_{k\in\mathcal{K}}\langle u,v_k\rangle x_k$ for any $u\in \mathcal{H}$ and the set $\{v_k\}_{k\in\mathcal{K}}$ is called a \textit{normalized tight frame.}
In particular, note that $c_1$ is merely a scaling factor, when the frame is tight.

While there are definitely similarities to the typical bases that one is used to, there are some important differences. 
Most importantly, $\mathcal{K}$ can be arbitrarily much larger than $\dim\mathcal{H}$, and $\{v_k\}_{k\in\mathcal{K}}$ could even contain repeated vectors. 
In particular, the following theorem, which is in essence a special case of an old result known as Naimark's dilation theorem \cite{czaja2008remarks}, has provided the authors with considerable insight. 
\begin{theorem}{(Na{i}mark)}
Every finite normalized tight frame $\{v_k\}_{k\in\mathcal{K}}$ for $\mathcal{H}$ is the orthogonal projection onto $\mathcal{H}$ of an orthonormal basis for a space of dimension $|\mathcal{K}|$, and vice versa.
\end{theorem}
A clear parallel to our case becomes evident from the following characterization. 
\begin{proposition}{\cite[Prop. 2.1]{waldron2018introduction}}\label{prop: tight frame NN*=zI}
A finite sequence $\{v_k\}_{k\in\mathcal{K}}$ in $\mathcal{H}$ is a tight frame for $\mathcal{H}$ with frame bound $c$ if and only if $NN^*=cI.$
Here, $N$ is the synthesis operator, i.e., $N(x)=\sum_{j=1}^n x_jv_j$ for $x\in \mathbb{C}^n$ and $N^*$ is the (dual) analysis operator, i.e. $N^*(x)_j=\langle x,v_j \rangle$. 
\end{proposition}

Indeed, we  observe a striking similarity to Proposition \ref{prop: NN*=kI}. 
We do not, however, immediately obtain equivalence. 
The main distinction between a general tight frame and the systems of lines that are of interest for this work, is that we require that the the frame vectors to have unit norm. 

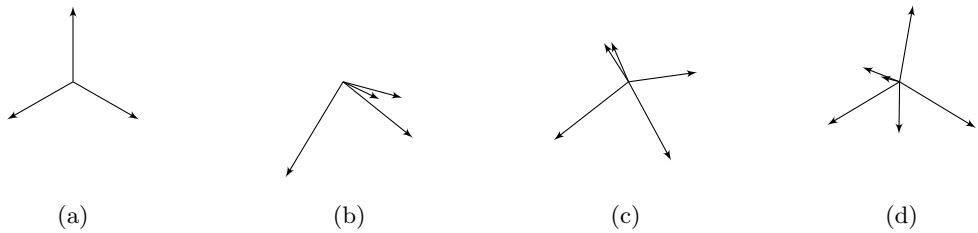
\begin{figure}[h!]
    \centering   
    \begin{subfigure}[b]{.24\textwidth}\centering
    \begin{tikzpicture}[scale=1.5]
\draw[arc] (0,0) -> (0,0.667);
\draw[arc] (0,0) -> (-0.5774,-0.333);
\draw[arc] (0,0) -> (0.5774,-0.333);
\node (e1) at (0,.7) {};
\node (e2) at (0,-.85) {};
    \end{tikzpicture}
    \caption{}
    \label{fig: mercedes frame}
    \end{subfigure}
    \begin{subfigure}[b]{.24\textwidth}\centering
    \begin{tikzpicture}[scale=1.5]
\draw[arc] (0,0) -> (-0.5038,-0.8436);
\draw[arc] (0,0) -> (0.6113,-0.4946);
\draw[arc] (0,0) -> (0.3207,-0.1543);
\draw[arc] (0,0) -> (0.5193,-0.1409);
\node (e1) at (0,.7) {};
\node (e2) at (0,-.85) {};
    \end{tikzpicture}
    \caption{}
    \end{subfigure}
    \begin{subfigure}[b]{.24\textwidth}\centering
    \begin{tikzpicture}[scale=1.5]
\draw[arc] (0,0) -> (-0.1521,0.3537);
\draw[arc] (0,0) -> (-0.6549,-0.5156);
\draw[arc] (0,0) -> (0.3735,-0.6956);
\draw[arc] (0,0) -> (0.6012,0.0840);
\draw[arc] (0,0) -> (-0.2169,0.3436);
\node (e1) at (0,.7) {};
\node (e2) at (0,-.85) {};

    \end{tikzpicture}
    \caption{}
    \end{subfigure}
    \begin{subfigure}[b]{.24\textwidth}\centering
    \begin{tikzpicture}[scale=1.5]
\draw[arc] (0,0) -> (-0.3261,0.1288);
\draw[arc] (0,0) -> (-0.0056,-0.4604);
\draw[arc] (0,0) -> (-0.1693,0.0461);
\draw[arc] (0,0) -> (0.6707,-0.4082);
\draw[arc] (0,0) -> (0.1173,0.6762);
\draw[arc] (0,0) -> (-0.6335,-0.3814);
\node (e1) at (0,.7) {};
\node (e2) at (0,-.85) {};
    \end{tikzpicture}
    \caption{}
    \end{subfigure}
    \caption{Three tight frames in $\mathbb{R}^2$ of order $3,4,5,6$, respectively. Only (a) is equal-norm.}
    \label{fig: tight frame strange}
\end{figure}

Note that indeed, this does not have do be the case, as is illustrated in Figure \ref{fig: tight frame strange}.
Thus, we are, in fact, looking for a special case of the above, which are called \textit{equal-norm} 
 tight frames.
 In particular, by Propositions \ref{prop: tight frame NN*=zI} and \ref{prop: NN*=kI}, it follows that $\{v_k\}_{k\in\mathcal{K}}$ is a unit-norm tight frame if and only if it corresponds to a two-eigenvalue gain graph via the usual construction.

The former case, in which the $\{v_k\}_{k\in\mathcal{K}}$ are an $\mathcal{A}_1$-set, is known as an \textit{equiangular} tight frame. 
These have been studied quite extensively, see e.g. \cite{szollHosi2014all, fickus2015tables, fickus2020mutually}, and will translate nicely to a two-eigenvalue gain graph on a complete graph. 
In addition to various sporadic examples, some infinite families of equiangular tight frames have been uncovered. 
For example, the vectors 
\begin{equation}
\label{eq: order 6 ETF}
    [0~~\tau~~\sigma]^\top,~~[\sigma~~0~~\tau]^\top,~~[\tau z~~\sigma~~0]^\top,~~[0~~\tau~~-\sigma]^\top,~~[-\sigma~~0~~\tau]^\top,~~[\tau z~~-\sigma~~0]^\top
\end{equation}
form an equiangular tight frame for arbitrary $z\in\T$ if $\tau=\sqrt{\frac{5+\sqrt{5}}{10}}$, $\sigma=\sqrt{\frac{5-\sqrt{5}}{10}}$. 
Interestingly, the corresponding gain graph is, in fact, a donut graph (see Definition \ref{def: donut}) with 
\[C=\begin{bmatrix}0&0&1\\1&0&0\\ 0&\overline{z}&0\end{bmatrix}. \]

Moreover, in his article \cite{renes2007equiangular}, Renes obtains an infinite family of equiangular tight frames of increasing order, based on the quadratic residues of prime powers of Gaussian primes (i.e., prime numbers congruent to $3\mod 4$). 
Below, we showcase an infinite family of gain graphs that may be distilled from this construction.

\begin{theorem}{\cite{renes2007equiangular}}\label{thm: renes family}
Let $p$ be a Gaussian prime and let $n:=p^z$, $z\in\mathbb{N}$. Then define the $n\times n$ matrix $M$ by 
\[M_{jh} = \begin{cases}0 & \text{ if } j=h \\ 1 &\text{ if } j\not=h \text{ and }~h-j \text{ is a quadratic residue in }GF(n) \\ -1 & \text{ otherwise.} \end{cases}\]
Now the matrix $A$, obtained from $M$ as
\[A = (n+1)^{-1/2}\left(I-J- i\sqrt{n}M \right)\]
has exactly two distinct eigenvalues.  
Specifically, it has an eigenvalue $\sqrt{n+1}$ with multiplicity $ (n-1)/2$ and an eigenvalue $-(n+1)^{3/2}/(n-1)$ with multiplicity $ (n+1)/2 $.
\end{theorem}
\begin{corollary}
For a given Gaussian prime $p$ and $n=p^z$, $z\in\mathbb{N}$, the $n\times n$ matrix $A$ as above characterizes a unit gain graph on $K_n$, with exactly two distinct eigenvalues. 
\end{corollary}

In case the frame vectors are allowed to be either orthogonal, or separated by angle $\alpha$, then much less previous work is readily available. 
Some work has been done on so-called $2$-angle tight frames or two-distance tight frames: see e.g., \cite{barg2015finite}.
However, their setting is more general (our case fixes one of the allowed angles to zero) and such works are generally concerned with constructions without orthogonal vectors. 

\subsection{SIC-POVM}
\label{subsec: sic povms}
As was briefly touched upon, an interesting link to two-eigenvalue gain graphs finds its origin in foundational quantum mechanics. 
Within the usual perimeters of quantum mechanics, it is commonly accepted that one can perform a number of manipulations to a system, but observing its state may compromise any further operation or destroy the state altogether. 
Moreover, since such observations are generally probabilistic in nature, one needs to be particularly careful when measuring the state of a quantum algorithm. 

In possession of certain defining qualities, the so-called \textit{symmetric, informationally complete, positive operator-valued measure} (SIC-POVM hereafter) is an interesting candidate to become the standard quantum measurement. 
Slightly paraphrasing \cite{geng2020minimal}, an IC-POVM is described by $m^2$ positive semi-definite operators $\{E_j\}_{j=1}^{m^2}$  that span the $m^2$-dimensional space of observables on an $m$-dimensional Hilbert space $\mathcal{H}$.
It is called an SIC-POVM if, in addition, it also satisfies the following three conditions:
\begin{enumerate}
    \item $E_j$ is rank one for all $j\in\{1,\ldots,m^2\}$, 
    \item Tr$E_jE_h = c$ for all $j\not=h$, $j,h\in\{1,\ldots,m^2\},$
    \item Tr$E_j = b$ for all $j\in\{1,\ldots,m^2\}$, 
\end{enumerate}
where $b$ and $c$ are nonnegative constants. 
In the below, we choose $b=1$ without loss of generality. 

Without too much effort, we see a clear relation of the above to an $\mathcal{A}_1$-set of lines that attains the bound in \eqref{eq: absolute bound}. 
Through that lens, the above is equivalent to a system of $m^2$ unit-norm vectors in an $m$-dimensional complex space, with pairwise equal inner products, in absolute value. 
Indeed, the $E_j$ are rank $1$ and positive semi-definite if they are the projectors $E_j = v_jv_j^*$, where Tr $v_jv_j^*=\|v_j\|^2=1$ for all $j\in[m^2]$ and Tr $v_jv_j^*v_hv_h^*=|v_j^*v_h|^2=c$  for all $j\not=h, j,h\in[m^2].$
Now, since the projectors $v_jv_j^*$ are a basis of $\mathbb{C}^{m\times m}$, (see below Proposition \ref{prop: sharper bound A2 set}) it follows that $c=(n-m)/(km)=1/(m+1)$ and that $\sum_{j=1}^{m^2}v_jv_j^*=mI,$ which by Proposition \ref{prop: NN*=kI} leads to a straightforward construction.

\begin{lemma}\label{lemma: tegg from sicpovm}
Let $\{E_1,E_2,\dots,E_n\}$ be a SIC-POVM in $\mathbb{C}^m$, so $n=m^2$. 
Let $N$ be the matrix whose columns are  $v_1,v_2,\dots,v_n$, where $E_j=v_jv_j^*$ for every $j\in[n]$.
Then $A=\sqrt{m+1}(N^*N-I)$ is the gain matrix of a two-eigenvalue gain graph.
\end{lemma}

The simplest example of a nontrivial SIC-POVM is obtained in $\mathbb{C}^2$ by the vectors that form the vectices of a regular tetrahedron in the Bloch sphere. 
Specifically, let 
\begin{equation}\label{eq: SICPOVM m=2}
    N_2 = \frac{1}{\sqrt{3}}\begin{bmatrix} \sqrt{3} & 1 & 1 & 1 \\ 
    0 & \sqrt{2} &\sqrt{2} \varphi & \sqrt{2} \varphi^2 \end{bmatrix},
\end{equation}
then $A(\Psi) = \sqrt{3}(N_2^*N_2-I) = W_4,$ which as we know has eigenvalues $\pm\sqrt{3}$. 

The first example that has not appeared in this work yet is obtained from a SIC-POVM of dimension $3$, is given by columns of the matrix $N_3$, below. 
\begin{equation}
    \label{eq: sic povm}
    N_3 = \frac{1}{\sqrt{2}}\begin{bmatrix}
    1 & 0 & \w & 1 & 0 & -1 & 1 & 0 & \ww \\
    \w & 1 & 0 & -1 & 1 & 0 & \ww & 1 & 0\\ 
    0 & \w & 1 & 0 & -1 & 1 & 0 & \ww & 1
    \end{bmatrix}
\end{equation}
\indent Given the results we have seen so far, the reverse construction is straightforward, though one would need a gain graph with very particular properties to work with.

\begin{proposition}
Let $\Psi=(K_n,\psi)$ be a gain graph of order $n=m^2$, for some $m\in \mathbb{N}$, whose spectrum is 
\[\Sigma_\Psi = \left\{(m-1)\sqrt{m+1}^{[m]},-\sqrt{m+1}^{[m^2-m]}\right\}.\] 
Then there exists an $N\in\mathbb{C}^{m\times n}$ such that $N^*N=I+\frac{1}{\sqrt{m+1}}A(\Psi);$ the columns of $N$ correspond to a SIC-POVM of dimension $m$.
\end{proposition}
As of yet, we do not posses the means to find gain graphs with the desired spectra of an order corresponding to a dimension for which existence of a SIC-POVM is an open problem. 
We feel that such an endeavor would be quite challenging, though further development of constructions such as the ones in Section \ref{sec: constructions} might lead to surprising results on this front.

\subsection{Mutually Unbiased Bases}
\label{subsec: mubs}
Another relevant concept from quantum information theory is that of \textit{Mutually Unbiased Bases} \cite{durt2010mutually}.
Where a SIC-POVM is effectively a maximum $\mathcal{A}_1$-set, a collection of MUBs is  a particular $\mathcal{A}_2$-set. 
Formally, 
two orthonormal bases $\{e_j\}_{j=1}^m$ and $\{f_h\}_{h=1}^m$ of $\mathbb{C}^m$  are said to be mutually unbiased if $|e_j^*f_h|^2=1/m$ 
for all $j,h\in[m].$ 

Much is known about MUBs. 
For instance, the maximum number of MUBs in $\mathbb{C}^m$ is $m+1$ when $m$ is a prime power, i.e., $m=p^z$, $z\in\mathbb{N}$, with $p$ prime. 
Yet, if $m$ is a different composite number, then the maximum number of MUBs is not known; even relatively small cases such as $m=6$ remain open\footnote{The general belief is that the maximum number of MUBs in  $\mathbb{C}^6$ is $3$.}. 

From the definition, one may already have observed the clear parallel to two-eigenvalue gain graphs. 
Indeed, if $N$ is the $ m \times (t\cdot m)$ matrix whose columns are $t$ MUBs, then  
$A=\sqrt{m}(N^*N-I)$
 defines a unit gain graph. 
Moreover, it follows from the definition that $NN^*=tI$, and thus Proposition \ref{prop: NN*=kI} applies, confirming that $A$ has exactly two distint eigenvalues $(t-1)\sqrt{m}$ and $-\sqrt{m}$.
Note that its underlying graph is the complete multipartite graph $K_{t\times m}$. 

We will briefly discuss some examples. 
For the smallest nontrivial case, the standard basis and the $4$ vectors 
\[\frac1{\sqrt{2}}\begin{bmatrix}1\\i^j\end{bmatrix}, j\in\{0,1,2,3\},\]
together form $3$ MUBs in $\mathbb{C}^2$. 
The corresponding two-eigenvalue gain graph, hereafter indicated by $K_{2,2,2}^{(\gamma)}$, is described by the following matrix:
\begin{equation}A(K_{2,2,2}^{(\gamma)}) = \begin{bmatrix}
0 & 0 &  -i                 & i & 1 & 1 \\
0 & 0 &   1                 & 1 & -1 & 1 \\
i & 1 &   0                 & 0 & -\bar{\gamma} & \gamma \\
-i & 1 &  0                 & 0 & -\gamma & \bar{\gamma} \\
1 & -1 & -\gamma & -\bar{\gamma}   & 0 & 0 \\
1 & 1 &   \bar{\gamma}            & \gamma & 0 & 0 
\end{bmatrix}.\label{eq: K222}\end{equation} 
In similar fashion, one may take consider $\mathbb{C}^3$. 
Here, the  vectors 
\begin{equation}\label{eq: MUB's in c3}
\begin{bmatrix}1&0&0 \end{bmatrix}^\top,\begin{bmatrix}0&1&0 \end{bmatrix}^\top,\begin{bmatrix}0&0&1 \end{bmatrix}^\top\text{~and~}
\frac1{\sqrt{3}}\begin{bmatrix}1 & \f^{j}& \f^h\end{bmatrix}^\top,~\text{where}~j,h\in\{0,1,2\}
,\end{equation}
form 4 MUBs in $\mathbb{C}^3$. 

It should be clear that when the set of MUBs does not attain the upper bound $m+1$ on its size, the above still holds true. 
For example, the toral tesselation graph $T_8^{(x)}$, where $x\in \T$, may alternatively be obtained from a pair of MUBs in $\mathbb{C}^4$. 
Specifically, if $N$ is the matrix whose columns are the standard basis 
appended with the basis formed by the vectors
\[
\frac12[1~1~1~-1]^\top, \frac12[1~1~-1~1]^\top, \frac12[1~-1~x~x]^\top, \frac12[-1~1~x~x]^\top
,\]
the graph characterized by $2(N^*N-I)$ is switching isomorphic to $T_8^{(x)}.$

As with the SIC-POVMs in the previous section, the relation between gain graphs with two specific eigenvalues and MUBs is an equivalence. 
\begin{proposition}
\label{prop: MUB}
Let $\Psi$ be a gain graph  on $K_{t\times m}$, with exactly two distinct eigenvalues $-\sqrt{m}$ and $(t-1)\sqrt{m}$. 
Then there exists an $N\in \mathbb{C}^{m\times tm}$ such that $A(\Psi)=\sqrt{m}(N^*N-I)$; the columns of $N$ form $t$ mutually unbiased bases in $\mathbb{C}^{m}$.
\end{proposition}

\subsection{Dismantling two-eigenvalue gain graphs}
\label{subsec: dismantling}
As was noted before, one must exercise considerable care when taking subgraphs of two-eigenvalue gain graphs. 
While it may at a first glance look intuitive to simply take a subset of the corresponding system of lines, the problem is that such a subset does not, in general, satisfy the necessary equation $NN^*=zI$.
However, a reliable way to obtain such subsystems is reminiscent of the dismantlability of certain association schemes \cite{martin2007imprimitive, lecompte2010equivalence, van2013uniformity}.

Suppose that the $m$-dimensional complex unit vectors $\{v_j\}_{j=1}^n$ correspond to a two-eigenvalue gain graph in the usual way. 
Then the matrix $N$, whose columns are the $v_j$, satisfies $NN^*=zI$ for some $z\in\mathbb{R}$. 
Now, if the columns of $N$ can be partitioned and concatenated into two matrices $N_1$ and $N_2$, where the former also constitutes a two-eigenvalue gain graph, then clearly $N_1N_1^*=yI$.
However, since $NN^* = \sum_{j=1}^n v_jv_j^* = N_1N_1^* + N_2N_2^*$, it follows that $N_2N_2^*=(z-y)I$, 
and thus, per the discussion following Proposition \ref{prop: NN*=kI}, the $N_2$ either corresponds to an empty graph, or to a two-eigenvalue gain graph. 

Clearly, the above may simply be repeated  
so long as one can find a subset of the vectors that satisfies the required equation. 
In particular, the $\{v_j\}_{j=1}^n$ may be partitioned into $s$ subsets that each satisfy $N_iN_i^*=z_iI$, with $z_i\in\mathbb{R}$ and $i=1,\ldots,s$.
Hence, every union of such subsets constructs a two-eigenvalue gain graph. 
Moreover, since every such gain graph is regular, the corresponding partition is an equitable one. 
The last conclusion, in particular, smells a lot like dismantlability.

Naturally, the hardest part is actually obtaining the desired subset of vectors, though the case of MUBs may serve as an illustrative example. 
Indeed, since MUBs may be partitioned into orthonormal bases, which of course each satisfy $N_iN_i^*=I$, every union of mutually unbiased bases corresponds to a two-eigenvalue gain graph. 
By the observation below Proposition \ref{prop: NN*=kI} and the discussion above, this holds in general, as long as the union of orthonormal bases forms an $\mathcal{A}_2$-set in $\mathbb{C}^m$. 

While these partitions are still not always easy to find, the above does yield an abundance of two-eigenvalue subgraphs of two-eigenvalue gain graphs that correspond to highly symmetric geometries, such as those discussed in Section \ref{subsec: examples lines}.

\subsection{A classification based on multiplicities}
\label{subsec: multiplicity lines}
Thus far, this section has mostly been concerned with drawing various parallels between the here considered gain graphs and various notions related to systems of lines in complex space. 
Existence of such systems, especially those whose cardinality is high compared to the dimension of the space in which they exist, is certainly no trivial issue.
In turn, said dimension corresponds to the multiplicity of the largest eigenvalue of our gain graphs.
As a consequence, we may use the discussed results to classify two-eigenvalue gain graphs with restricted multiplicity, by considering the systems of lines in low-dimension spaces.

The following interesting fact, known as the Cvetkovi{\'c} bound \cite[Thm.~3.5.1]{brouwer2011spectra}, is a particularly useful tool in the approach that is taken below. 
\begin{lemma}
Let $\Psi$ be a gain graph with spectrum $\left\{\theta_1^{[m]},\theta_2^{[n-m]}\right\}$. 
Then the largest coclique in $\Psi$ has size at most $m$. 
\end{lemma}
\begin{proof}
Suppose $\Psi$ contains a coclique of size $m+1$. Then, by eigenvalue interlacing, $\Psi$ has $\lambda_{m+1}\geq 0,$ which is a contradiction. 
\end{proof}
Recall that least multiplicity $1$ occurs only for gain graphs that are switching equivalent to a complete graph, so we proceed to the smallest interesting case.
\subsubsection{Multiplicity $2$}
\label{subsec: k=2}
Let us classify the gain graphs $\Psi$ whose spectrum is exactly 
\begin{equation}
    \left\{\theta_1^{[m]}, \theta_2^{[n-m]}\right\}, \text{~where~} m=2. 
\end{equation}
As noted above, this corresponds to a system of lines in $\mathbb{C}^2$. 
Applying the bound in \eqref{eq: absolute bound}, we obtain that $n\leq 4$ if the gain graph is complete and $n\leq 6$ otherwise. 
In the former case, one trivially obtains $K_3$; the case $n=4$ is exactly the unique SIC-POVM in \eqref{eq: SICPOVM m=2}, whose corresponding gain graph is\footnote{Note that any unitary transformation of the system of lines does not change the corresponding gain graph. That is, if $U$ is a unitary matrix, then $M:=UN$ and $N$ represent the same gain graph because $M^*M=N^*N$.}  $W_4$.  

Next, consider the case that $\Gamma(\Psi)$ is not complete; i.e., $\Psi$ corresponds to an $\mathcal{A}_2$-set.
Without loss of generality, we choose the first vector to be the standard unit vector $e_1$. 
Then, since $\Psi$ is not complete, at least one vector must be orthogonal to $e_1$, so (without loss of generality) we take the second unit vector $e_2$. 
Since $m=2$, there is no vector that is orthogonal to both $e_1$ and $e_2$, so it follows that $\Psi$ has (constant) degree $k=n-2$. 
This, in turn, means that $\Psi$ is complete multipartite $K_{t\times 2}$, 
and thus (by Proposition \ref{prop: MUB}) $\Psi$ corresponds to a pair ($n=4$) or a set of three ($n=6$) MUBs. 
This yields precisely a gain graph that is switching isomorphic to $IG(W_2)$ or to $K_{2,2,2}^{(\gamma)}$, respectively.

\subsubsection{Multiplicity $3$}
\label{subsec: k=3}
We may apply the same line of questioning for graphs with spectrum $\{\theta_1^{[3]},\theta_2^{[n-3]}\}.$
If $k=n-1$, then the corresponding line system is an equiangular frame, and the absolute bound $n\leq 9$ applies. 
These have been classified in dimension $3$ by \citet{szollHosi2014all}, and the results are summarized in Table \ref{table: szollosi results}.
\begin{table}[h]
    \centering
    \begin{tabular}{rl}
    \toprule
      $n$ &  \\ \cmidrule{1-2}
      3   & Standard basis \\
      4   & Regular simplex  \\ 
      5   & - \\
      6   & Corresponds to \eqref{eq: order 6 ETF}\\
      7   & Corresponding graph obtained from Theorem \ref{thm: renes family}\\
      8   & - \\
      9   & Corresponds to  \eqref{eq: sic povm}\\
         \bottomrule
    \end{tabular}
    \caption{Classification $\mathcal{A}_1$-sets in $\mathbb{C}^3$ from \citet{szollHosi2014all}, up to equivalence.}
    \label{table: szollosi results}
\end{table}

In case $k=n-2$, then $\Gamma(\Psi)=K_{t\times 2}$, where $t\geq 3$ since it was assumed throughout that $m\leq n/2$.  
In other words, the vectors come in orthogonal pairs, with vectors from distinct pairs being separated by some angle $\alpha$.
Note that for given $n,m,$ and $k$, we may compute the value of $\theta_2$, which, as discussed before, determines $\alpha$.
Furthermore, note that since $K_{t\times 2}$ does not have an eigenvalue $-3(2t-2)/(2t-3),$ Proposition \ref{prop: sharper bound A2 set} implies that $n\leq 9.$ 

Using the above, we are left with just two cases: either $n=6$ or $n=8$. 
If $n=6$, then $k=4$ and $m=3$, which means that $\alpha=1/2$.
Again, take the first two standard unit basis vectors as a starting point. 
Then any further candidate $v$ is of the form $v = \begin{bmatrix}
x/2 & y/2 & z/\sqrt{2}
\end{bmatrix}^\top$, for $x,y,z\in\T$. 
It is not hard to see that (up to equivalence) the system  represented by the vectors
\begin{equation}
    \label{eq: mult3MUB}
    \begin{bmatrix}1 \\0  \\ 0 \end{bmatrix},
            \begin{bmatrix}0 \\1  \\ 0 \end{bmatrix},
            \begin{bmatrix}1/2 \\1/2  \\ 1/\sqrt{2} \end{bmatrix},
            \begin{bmatrix}1/2 \\1/2   \\-1/\sqrt{2} \end{bmatrix},
            \begin{bmatrix}1/2 \\-1/2  \\ z/\sqrt{2} \end{bmatrix},
            \begin{bmatrix}1/2 \\-1/2  \\ -z/\sqrt{2} \end{bmatrix}
\end{equation}
satisfies the requirements for any $z\in\T$. 
The corresponding two-eigenvalue gain graph is, in fact, the order-$6$ toral tesselation graph $T_6^{(z)}$. (See the end of Section \ref{sec: constructions}.)

What is left is the the case $n=8$. 
We will formally show that $(n,m,k)=(8,3,6)$ does not yield any two-eigenvalue gain graphs.
\begin{lemma}\label{lemma: n=8 m=3 k=6 no go}
There is no gain graph $\Psi$ with spectrum $\{\theta_1^{[3]},\theta_2^{[5]}\}$ and $k=6$.
\end{lemma}
\begin{proof}
Since $(n,m,k)=(8,3,6)$ it follows that $\alpha=\sqrt{5/18}$ and $\Gamma(\Psi)=K_{4\times 2}.$
Hence, the corresponding line system contains four pairs $(v_j,v_{j+1})$, $j\in\{1,3,5,7\}$,  of orthogonal lines, such that any two lines from distinct pairs are separated by angle $\alpha$. 
Without loss of generality, choose the first pair to be $(v_1,v_2)=(e_1,e_2)$. 
Then any other candidate is of the form 
\[v_j=\begin{bmatrix}\alpha x_j & \alpha y_j & \beta z_j\end{bmatrix}^\top \text{ with } x_j,y_j,z_j\in\T \text{ for } j=3,\ldots,8 \text{ and where } \beta=\frac{2}{3}.\]
As before we may assume w.l.o.g. that $x_j=1$ for all $j=3,\ldots,8$.
Now, if $N$ is the matrix whose columns are $v_1,\ldots,v_8$ and $\Psi$ is a two-eigenvalue gain graph, then $NN^*$ is a multiple of the identity. 
In other words, the rows of $N$ are orthogonal. 
It follows straightforwardly that 
\begin{equation}
    \label{newproof410 1}
    \sum_{j=3}^8y_j=0 \text{  ~and~  }\sum_{j=3}^8z_j=0. 
\end{equation}
Furthermore, let $y:=\bar{y}_3y_4$ and $z:=\bar{z}_3z_4$. 
Then $v_3^*v_4=\alpha^2+\alpha^2y+\beta^2z=0$ if and only if $y=-1-\beta^2\alpha^{-2}z.$
Since $|y|=|z|=1$, $z$ is on the intersection of two real-centered circles on the complex plane, and thus there is (at most) one conjugate pair of solutions for $z$. 
Suppose that $\tilde{z}$ is such a solution, and set $\tilde{y}=-1-\beta^2\alpha^{-2}\tilde{z}$.
By symmetry, we also have\footnote{Up to conjugation, though since $\overline{\bar{y}_5y_6}=\bar{y}_6y_5$, equality is assumed without loss of generality.} $\bar{y}_5y_6 = \bar{y}_7y_8=\tilde{y}$ and $\bar{z}_5z_6 = \bar{z}_7z_8=\tilde{z}$. 
Plugging this into \eqref{newproof410 1} yields $\sum_{j=3}^8y_j=0$ if and only if $y_3+y_5+y_7=0$, and similarly $z_3+z_5+z_7=0$. 
Here, it is used that $\tilde{y}\not=-1$, since this would contradict $|z|=1$. 

Finally, w.l.o.g. assume $y_3=z_3=1$. 
Then the equations $|y_5|=|y_7|=1$ and $y_5+y_7=-1$ must simultaneously hold, which implies $y_5=\pm\varphi$ and $y_7=\bar{y}_5$; similarly $z_5=\pm\varphi$ and $z_7=\bar{z}_5$. 
But now $|v_3^*v_5|=|\alpha^2 + \alpha^2y_5 + \beta^2z_5|\not=\alpha$, which is a contradiction. 
\end{proof}

Finally, in case $k\leq n-3$,  we find that $\Psi$ must contain a coclique of order (at least) 3. 
\begin{lemma}
\label{lemma: mult 3 coclique}
Let $\Psi$ be an order-$n$ gain graph with constant degree $k\leq n-3$, and spectrum $\left\{\theta_1^{[3]},\theta_2^{[n-3]}\right\}$.
Then $\Psi$ contains a coclique of order $3$.
\end{lemma}
\begin{proof}
We will be reasoning with the system of lines that corresponds to $\Psi$ in the usual way.
Suppose that $\Psi$ does not contain a coclique of order $3$ and let $e_1$ be the first standard unit vector. 
If $k\leq n-3$, then there are (at least) two vectors that are orthogonal to $e_1$. 
Moreover, said vectors may not be orthogonal to one another, since this would imply existence of an order-$3$ coclique. 
Hence, they must be separated by an angle 
$\alpha = \sqrt{(n-3)/(3k)}$.
Without loss of generality, choose the second unit vector $e_2$ and $v=[0~\alpha~z],$ where $|z|=\sqrt{1-\alpha^2}.$

Now, we may repeat the argument, since, in addition to $e_1$, there is at least one more vector that is orthogonal to $e_2$. 
Furthermore, said vector must make an angle $\alpha$ with $e_1$.
Without loss of generality, choose $w = [\alpha ~0 ~z']$, where $|z'|=\sqrt{1-\alpha^2}.$
However, now $|v^*w|=|z^*z'|=1-\alpha^2$, which implies $\alpha = 1-\alpha^2$, and thus $\alpha=(\sqrt{5}-1)/2$, which is a contradiction for all possible $(n,k)$.
\end{proof}
Essentially, by using the symmetry in the argument above, it follows that any vertex in a two-eigenvalue gain graph $\Psi$ with the desired spectrum is contained in an order-$3$ coclique. 
Then, it is not too hard to see that the inclusion of such an order-3 coclique implies that the corresponding line system consist of MUBs.
\begin{lemma}
\label{lemma: mult 3 MUBs}
Let $\Psi$ be a two-eigenvalue gain graph with least multiplicity $3$. 
If $\Psi$ contains an order-3 coclique then the usual corresponding line system consists of mutually unbiased bases.
\end{lemma}
\begin{proof}
Consider the corresponding line system. 
Since $\Psi$ contains a coclique of size $3$, the system contains an orthonormal basis; without loss of generality, assume that it is the standard basis. 

Now, suppose another vector $v$ has first entry $0$. 
By considering $|e_2^* v|$ and $|e_3^* v|$, which should both equal $\alpha$, it follows that $\alpha=\sqrt{1/2}.$
Moreover, it also follows that $k<n-3$, since there are 3 vectors that are orthogonal to $e_1$. (Specifically, $e_2,e_3$ and $v$.)
Hence, there must be a third vector with second entry $0$, which (by repeating the above) must satisfy $|e_1^* w|=|e_3^* w|=\alpha$. 
But now $|v^* w|=|(e_3^*v)(e_3^*w)|=1/2\not=\alpha$, which is a contradiction.
Hence, no vectors other than the standard basis may have zero entries, and thus $k=n-3$. 
The conclusion now follows since every vertex is contained in a unique 
order-$3$ coclique.
\end{proof}

It follows immediately that the remaining gain graphs with the desired spectrum have order $n=6,9$ and $12$, and may be obtained from the collection of MUBs described in \eqref{eq: MUB's in c3}.  

To conclude our classification, we include a summarizing Theorem.
Named graphs will be listed as such; unnamed examples are referenced using the equation that contains it or an equivalent system of lines. 
\begin{theorem}
\label{thm: mult at most 3}
All two-eigenvalue gain graphs with least multiplicity at most $3$ are switching isomorphic to one of the gain graphs in Table \ref{tab: mult summary}.
\end{theorem}
\begin{table}[h!]
\centering
\begin{tabular}{@{}clllccllllc@{}}
\toprule
$m$ & \multicolumn{1}{l}{Graph} & \multicolumn{1}{l}{Order} & \multicolumn{1}{l}{$k$} &\multicolumn{1}{l}{DS} & &\multicolumn{1}{l}{Graph} & \multicolumn{1}{l}{Order} & \multicolumn{1}{l}{$k$} &\multicolumn{1}{l}{DS} \\ \cmidrule{2-5}\cmidrule{7-10}
 \multirow{1}{*}{$1$}                   & $K_n$                    &  $n$                       & $n-1$   & * &
                                        &             &                     &             &  \\ \\
 \multirow{2}{*}{$2$}                   & $IG(W_2)$                &  $4$                     & $2$     & * &
                                        & $K_{2,2,2}^{(\gamma)}$                  &  $6$                     & $4$         & * \\
                                        & $W_4$                        &  $4$                     & $3$ & * 
                                        &             &                     &           &    \\
                     &                           &                           &            
                                        &             &                     &        &       \\               
 \multirow{4}{*}{$3$}                   & 2 MUBs from \eqref{eq: MUB's in c3}&  $6$                     & $3$               &   *   & 
                                        &   3 MUBs from \eqref{eq: MUB's in c3}               &  $9$      & $6$     &       *        \\ 
                                        & $T_6^{(x)}$                   &  $6$                     & $4$                &     & 
                                        & Equation \eqref{eq: sic povm} &  $9$                     & $8$       &        *        \\ 
                                        & Equation \eqref{eq: order 6 ETF}  &  $6$              & $5$                &     &  
                                        & 4 MUBs from \eqref{eq: MUB's in c3} &  $12$              & $9$       &     *         \\ 
                                        & Theorem 4.4                          &  $7$             & $6$                      & * 
                                        &                                     &                   &         &        \\\bottomrule
\end{tabular}
\caption{Classification of all two-eigenvalue gain graphs with least multiplicity at most $3$. A star in the DS column indicates that any connected, cospectral gain graph is switching isomorphic. }
\label{tab: mult summary}
\end{table}

Naturally, one could choose to increase the multiplicity further and apply more or less the same arguments again; the final series of arguments (concerning $k\leq n-m$) in particular appears as if it would carry over with little to no issues. 
However, since we feel that such a discussion would provide little new insight, we choose to move on.

\subsection{Other low dimension examples}
\label{subsec: examples lines}
We conclude this section by showcasing a number of two-eigenvalue gain graphs that arise from various other well-studied combinatorial objects.
The constructions are similar to those applied earlier in this section, though their dimensions and the corresponding graph orders are higher.
\subsubsection{The Witting polytope}
Real polytopes have been generalized to complex Hilbert spaces for quite some time. 
While precise definitions do not exist for the general case, the regular complex polytopes have been completely characterized by Coxeter. 
We forego the details, though it would be fair to say that these geometries are highly symmetric, which enables us to translate (parts of) them to the desired systems of lines.

Consider the Witting polytope \cite{coxeter1991regular} in $\mathbb{C}^4$. 
Its 240 vertices occur in $40$ $1$-dimensional subspaces, which form an $\mathcal{A}_2$-set meeting the absolute bound. 
In particular, take the 4 standard basis vectors along with
\begin{align*}
     \frac1{\sqrt{3}}\begin{bmatrix}1~ 0~ -\f^{j}~ -\f^{h}\end{bmatrix}^\top, 
     \frac1{\sqrt{3}}\begin{bmatrix}1~ -\f^{j}~ 0~ \f^{h}\end{bmatrix}^\top,
     \frac1{\sqrt{3}}\begin{bmatrix}1~ \f^{j}~ \f^{h}~ 0\end{bmatrix}^\top,\text{~and~} \\
     \frac1{\sqrt{3}}\begin{bmatrix}0~ 1~ -\f^{j}~ \f^{h}\end{bmatrix}^\top, \text{ with } j,h\in\{0,1,2\}.
\end{align*}
Then the matrix $N$, whose columns are the vectors above, satisfies $NN^*=10I,$ so that Proposition \ref{prop: NN*=kI} applies.
Indeed, $A=\sqrt{3}(N^*N-I)$ characterizes an order-$40$ unit gain graph whose spectrum is $\left\{9\sqrt{3}^{[4]},-\sqrt{3}^{[36]}\right\}$, and whose underlying graph is the complement of the symplectic generalized quadrangle of order $3$. 

Finally, note that the $40$ vectors above may be partitioned into ten orthonormal bases, which form a spread in said quadrangle. 
Following the discussion in Section \ref{subsec: dismantling}, one may use this partition to get two-eigenvalue gain graphs with 
spectrum $\left\{(t-1)\sqrt{3}^{[4]},-\sqrt{3}^{[4(t-1)]}\right\}$, $t\in\{2,\ldots,10\}$.

\subsubsection{A rank-$5$ complex reflection group}
In the same vein, one may draw on the collection of complex reflection groups and distill a unit gain graph from its hyperplanes of symmetry. 
In particular, we consider the group named ST33 in the (complete) classification by Shephard and Todd \cite{shephard1954finite}.

The desired collection of vectors 
may roughly be divided into two parts.
The former part consists of all vectors obtained from 
\[\frac{1}{\sqrt{2}}[1~~-\varphi^j~~0~~0~~0]^\top, j=0,1,2,\]
by permuting the first four entries in all possible ways, such that the leftmost nonzero entry is\footnote{Note that the vector obtained by interchanging the two nonzero entries is a member of the same hyperplane.} $1/\sqrt{2};$ this yields $18$ vectors belonging to distinct hyperplanes in $\mathbb{C}^4$.
The second part, containing the remaining $27$ vectors, is described by 
\[\frac1{\sqrt{6}}[1~~ \f^{j_1}~~ \f^{j_2}~~ \f^{j_3}~~ \sqrt{2}\f^{-j_1-j_2-j_3}]^\top, \text{~where~} j_1,j_2,j_3\in\{1.2.3\}. \]
This yields a total of 45 unit vectors, whose pairwise inner products have absolute value either $0$ or $\frac{1}{2}$, such that every vector is orthogonal to  exactly 12 others. 
Thus, the matrix $A=2(NN^*-I)$ characterizes a unit gain graph with spectrum $\left\{16^{[5]},-2^{[40]}\right\}$. 

Note that the underlying graph is the generalized quadrangle of order $(4,2)$, which does not have a spread. 
However, we have found a partition of the vector set that contains six orthonormal bases; this corresponds to a partial spread.
By considering unions of such bases and their complementary sets, we obtain induced subgraphs whose spectra are $\left\{2(t-1)^{[5]},-2^{[5(t-1)]}\right\}$ for $t=2,\dots,8$.
Interestingly, the corresponding gain matrix has entries in $\T_6$, so it (and its induced subgraphs) could be interpreted as a signed digraph. 

\subsubsection{The Coxeter-Todd lattice}
Finally, we consider the famous Coxeter-Todd lattice in $\mathbb{C}^6$, which finds its origin in the hexacode, discussed at the start of this section. 
It has various equivalent descriptions \cite{conway1983coxeter}, which give rise to different interesting two-eigenvalue gain (sub)graphs. 
In particular, these graphs attain the absolute bound \eqref{eq: absolute bound} in terms of order $n$ with respect to the multiplicity $6$.
The descriptions may be distinguished by their base\footnote{The lattice is said to be represented in the $b$-base if all absolute values of the pairwise inner products of the coordinate vectors are divisible by $b$, before scaling \cite{conway1983coxeter}.}. 

In the 2-base, take all 15 projectively distinct hexacodewords of weight $4$, and take all of their variations by multiplying at most 3 nonzero entries of every such codeword by $-1$. 
Scaling them down yields a collection of 120 distinct unit vectors, that may be appended with a standard unit basis of $\mathbb{C}^6$ to find 126 unit vectors whose pairwise inner products have absolute value either $0$ or $\frac{1}{2}$.
Thus, as before, we find a two-eigenvalue gain graph. 
Note, specifically, that the before-mentioned example on 15 vertices clearly occurs as an induced subgraph of this construction, as its corresponding system of lines is a subset of the 126 vectors required here. 

In fact, the 126 vectors decompose into 21 orthogonal bases.
This is easy to see, since the gain graph above has an order-$126$ underlying graph that is the complement of the strongly regular graph that appears in \citet{brouwer2021fragments} as $NO^-_6(3)$.
According to \cite{brouwer2021fragments}, the complement of $NO_6^-(3)$ has chromatic number $21$, which implies that the desired decomposition exists.
This, as before, yields two-eigenvalue gain subgraphs of order $6t$, for $t=2,\ldots,21.$

One might also consider the 3-base parallel. 
Indeed, take the 45 vectors obtained from 
\[[i\sqrt{3}~-i\sqrt{3}\f^c~0~0~0~0]^\top, c\in[3]\]
by permuting its entries in such a way that the first nonzero entry is strictly imaginary, and append with the 81 vectors 
\[[1~\f^{j_1}~\f^{j_2}~\f^{j_3}~\f^{j_4}~\f^{-j_1-j_2-j_3-j_4}]^\top\text{ with } j_1,\ldots,j_4\in[3].\]
It is easily verified that the pairwise inner products of these 126 vectors have absolute value either $0$ or $3$, and thus the usual construction applies, after scaling.

It should be noted that the collections of 45 and 81 vectors also construct two-eigenvalue gain graphs (with two-eigenvalue subgraphs), that are notably different from the subgraphs of order $6t$, above. 
In fact, the 81 vectors of weight 6 have yet another interesting link to other  combinatorial objects constructed by \citet{vanlint1981construction}, such as their partial geometry. 
There is a clear correspondence of the above $81$ vectors and the dual code of \cite[Construction 2]{vanlint1981construction}. 
The different inner products (up to conjugation) of our vectors correspond to the weights in this dual code \cite[Table III]{vanlint1981construction}, and define a $4$-class fusion scheme of the $8$-class cyclotomic association scheme on $GF(81)$, which can be further fused to an amorphic association scheme \cite{vandam2010some}. 
We should note that  \citet{roy2014complex} have obtained many results such as the above, where the inner products give rise to various association schemes. 
These most interesting constructions are called spherical $t$-designs, which in a sense generalize the above.

Lastly, we draw from the 4-base variant. 
Take the 96 distinct vectors obtained from
\[[i\sqrt{3}~(-1)^{j_1}~(-1)^{j_2}~(-1)^{j_3}~(-1)^{j_4}~(-1)^{-j_1-j_2-j_3-j_4}]^\top, \text{ with }j_1,\ldots,j_4\in[2],\]
by permuting its entries and append with the 30 obtained as the pairwise linearly independent permutations of 
\[[2~\pm2~0~0~0~0]^\top.\]
Then the pairwise inner products have absolute value $0$ or $4$, and the usual construction applies.

It turns out that each of these constructions yields a gain graph that belongs to the same switching equivalence class. 
Will will not offer formal argumentation, but it is easily verified by computer. 
Note, moreover, that this is quite unsurprising, since the lines were drawn from various descriptions of the same group. 
Additionally, it turns out that each of the obtained gain graphs once again has all of its nonzero entries in $\T_6$, thus admitting a signed digraph interpretation.

\section{Two eigenvalues and small degree}          
\label{sec: small degree}
Having classified all two-eigenvalue gain graphs with bounded multiplicity, a relatively small collection of admissible graphs was obtained. 
We will now take a different perspective in bounding the degree of a candidate gain graph to a small number; this will warrant a combinatorial approach in which we will systematically investigate the potential underlying graphs and the corresponding gain functions that may act on their edges. 
Interestingly, for all degrees $k\leq 4$, we obtain implicit bounds on the order of our candidates. 

We should note that this classification in the context of Hermitian adjacency matrices and Eisenstein matrices has essentially already been done by \citet{greaves2012cyclotomic}. 
Foregoing the complete graphs, his constructions all have $a=0$. 
However, as might be expected by now, we cannot generally make this assumption for gain graphs. 
Indeed, Section \ref{sec: systems of lines} contains many examples to the contrary, such as $K_{2,2,2}^{(\gamma)}$, 
which was obtained from $3$ mutually unbiased bases in $\mathbb{C}^2$ and  has distinct eigenvalues $2\sqrt{2}$ and $-\sqrt{2}$.

A particularly useful insight that will be helpful in the classification below, is the following. 
\begin{lemma}\label{lemma2cn} Let $D$ be a connected two-eigenvalue gain graph. Then any two vertices at mutual distance $2$ have at least two common neighbors.
\end{lemma}
\begin{proof} If there would be exactly one common neighbor $w$ between non-adjacent vertices $u$ and $v$, then $(A^2)_{uv} =A_{uw}A_{wv} \neq 0$. This contradicts the equation $A^2=aA+kI$.
\end{proof}
Additionally, by the following observation, we may substantially limit the decision space in case $\Psi$ is triangle-free. 
\begin{lemma}
\label{lemma: triangle-free implies a=0}
Let $\Psi$ be a connected two-eigenvalue gain graph. 
If $\Psi$ is triangle-free, then $a=0$ and thus $\theta_1=-\theta_2$. 
\end{lemma}
\begin{proof}
$a\not=0$ implies that $(A^2)_{uv}\not=0$ for connected vertices $u,v$, and thus there is a walk $u\to w\to v$. But then $\Psi[\{u,w,v\}]$ is a triangle.  
\end{proof}

Before we get into the actual classification, we would like to recall that it is assumed throughout that $a=\theta_1-\theta_2\geq 0$; that is, the eigenvalue whose multiplicity is lower is assumed to be positive. 
While this is a restrictive assumption, each of the graphs that are subsequently excluded may be obtained by multiplying one of the obtained graphs by $-1$.
Hence, nothing is effectively lost.

\subsection{Degree $2$}
The first relevant case to consider is, of course, $k=2$. 
Without much effort, we show that there are exactly two switching equivalence classes that admit to the imposed requirements.

\begin{theorem}
\label{thm: degree two}
Let $D$ be a connected unit gain graph with degree $k=2$ that has two distinct eigenvalues. Then $D$ is switching isomorphic to $K_3$ or  $\IG(W_2)$.
\end{theorem}
\begin{proof}
Suppose that $D$ is not a balanced triangle. 
Then both eigenvalues of $D$ have multiplicity at least $2$, and thus $n\geq 4$. 
Now, since $k=2$, then by Lemma \ref{lemma2cn} we have $\Gamma(D)\cong C_4$. 
By Lemma \ref{lemma: make a tree equal}, $3$ out of $4$ edges may be assumed to have gain one; setting $A^2=2I$ then easily yields the desired conclusion. 
\end{proof}

\subsection{Degree $3$}
\begin{figure}[t]
    \centering
    \begin{subfigure}[b]{.35\textwidth}\centering
    \begin{tikzpicture}[scale=1.1]
        \def\x{.6}
        \def\r{2}
            \node[vertex] (1) at (\x,\x) {};
            \node[vertex] (2) at (\x,-\x) {};
            \node[vertex] (3) at (-\x,-\x) {};
            \node[vertex] (4) at (-\x,\x) {};
            
            \node[vertex] (5) at (\r*\x,\r*\x) {};
            \node[vertex] (6) at (\r*\x,-\r*\x) {};
            \node[vertex] (7) at (-\r*\x,-\r*\x) {};
            \node[vertex] (8) at (-\r*\x,\r*\x) {};
            
            \draw[tree] (1) to node{} (2);
            \draw[tree] (3) to node{} (2);
            \draw[tree] (1) to node{} (4);
            \draw[negedge] (3) to node{} (4);
            
            \draw[negedge] (5) to node{} (6);
            \draw[negedge] (7) to node{} (6);
            \draw[negedge] (5) to node{} (8);
            \draw[edge] (7) to node{} (8);
            
            \draw[tree] (1) to node{} (5);
            \draw[tree] (2) to node{} (6);
            \draw[tree] (3) to node{} (7);
            \draw[tree] (4) to node{} (8);
    \end{tikzpicture}
    \caption{$\ND(\IG(W_2))$}
    \label{fig: signed cube}
    \end{subfigure}
    \begin{subfigure}[b]{.35\textwidth}\centering
    \begin{tikzpicture}[scale=1.35]
        \def\x{.6}
        \def\r{2}
        
            \node[vertex] (1) at (0.5*1.7321*\x,0.5*\x) {$5$};
            \node[vertex] (2) at (-0.5*1.7321*\x,0.5*\x) {$6$};
            \node[vertex] (3) at (0,-\x) {$4$};
            \node[vertex] (4) at (0.5*1.7321*\x*\r,-0.5*\x*\r) {$3$};
            \node[vertex] (5) at (-0.5*1.7321*\x*\r,-0.5*\x*\r) {$2$};
            \node[vertex] (6) at (0,\x*\r) {$1$};
            
            \draw[negarc] (1) to node{} (4);
            \draw[negarc] (5) to node{} (1);
            \draw[tree] (1) to node{} (6);
            
            \draw[negarc] (4) to node{} (2);
            \draw[negarc] (2) to node{} (5);
            \draw[tree] (2) to node{} (6);
            
            \draw[tree] (3) to node{} (4);
            \draw[tree] (3) to node{} (5);
            \draw[tree] (3) to node{} (6);

    \end{tikzpicture}
    \caption{$IG(W_3)$}
    \label{fig: k33}
    \end{subfigure}
    \caption{Illustrations for Theorem \ref{thm: degree 3}. A dashed edge has gain $-1$; a dashed arc has gain $\varphi$.}
    \label{fig: degree 3 illustrations}
\end{figure}
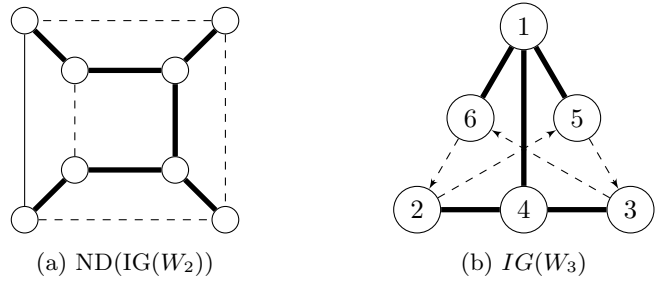

Increasing the degree of the considered candidates to $3$ gives us some more freedom, though the collection of switching equivalence classes is still limited to $4$. 
The desired classification is obtained with relative ease, by application of the process described above.
As was announced earlier, we find that the examples in Proposition \ref{prop: deg 3 2evggs} form a complete list. 
\begin{theorem} 
\label{thm: degree 3}
Let $\Psi$ be a connected unit gain graph with degree $k=3$ that has two distinct eigenvalues. Then $D$ is switching isomorphic to $K_4$ or one of the graphs $W_4$, $\IG(W_3)$, or the signed $3$-cube $\ND(\IG(W_2))$.
\end{theorem}
\begin{proof} 
Suppose that $\Psi$ is not switching equivalent to $K_4$ and let $A:=A(\Psi)$. 
We distinguish two cases: either $\Psi$ contains a triangle, or  it does not. 

    Suppose that $\Psi$ contains a triangle. Then without loss of generality, we may assume that $A_{12}=A_{13}=A_{14}=1$ and $A_{23}\neq 0$. 
    Then, equating the first column of $A^2$ to the corresponding entries of $aA+kI$, we find
    \[
    A_{23} + A_{24} = a ~\wedge~
    \overline{A_{23}}+A_{34}=a ~\wedge~
    \overline{A_{24}}+\overline{A_{34}}=a.
    \]
    Now, since $a$ is real, it follows that $A_{23}=\overline{A_{24}}=A_{34}$ and thus that $\Gamma(\Psi)=K_4$.
    Indeed, since $\Psi$ is connected and the first four vertices all have degree $3$, no further vertices can be added. 
    Moreover, since $\Psi$ is not switching isomorphic to $K_4$, its eigenvalues must be $\pm \sqrt{3}$ and thus $a=0$. 
    Finally, using the equations above, this implies $A_{23}=\pm i$; both choices yield $W_4$.  

    If $\Psi$ does not contain triangles, then by Lemma \ref{lemma: triangle-free implies a=0}, $a=0$. 
    Moreover, by Lemma \ref{lemma2cn}, every two vertices at distance 2 have at least two common neighbors. 
    There are now two subcases to distinguish, and each leads to one graph.
First, suppose that all pairs of vertices at distance $2$ have precisely two common neighbors.
Then $\Gamma(\Psi)$ is the cube \cite{brouwer2006classification}, of which the edges marked fat in Figure \ref{fig: signed cube} may be fixed to gain $1$. 
Now, note that every two non-adjacent vertices of every face are connected with exactly two walks of length $2$. 
Since the off-diagonal entries of $A^2$, which must be zero, are given by the sum of the gains of such a pair of walks, the gains of the non-fixed edges are all determined by the equation $A^2=3I.$
The resulting signed graph corresponds to $\ND(\IG(W_2))$, the signed cube.

In the other case, there must be two non-adjacent vertices, say $1$ and $2$, that share all three neighbors.
It follows that $\Gamma(\Psi)\cong K_{3,3}$; label such that $1,2,3$ are pairwise nonadjacent.  
Without loss of generality, we set $A_{14}=A_{15}=A_{16}=A_{42}=A_{43}=1.$ 
Then, since we have \[(A^2)_{12}=\phi(1\to4\to2)+ \phi(1\to5\to2)+\phi(1\to6\to2) = 1+A_{52}+A_{62}=0,\]
it follows that $A_{52}=\overline{A_{62}}$ and either $A_{52}=\f$ or $A_{52}=\ff$.
Moreover, since the same argument holds when $(A^2)_{13}$ is considered, we also have $A_{53}=\overline{A_{63}}$ and either $A_{53}=\f$ or $A_{53}=\ff$.
Finally, since $(A^2)_{23}=0$, we obtain $A_{52}=\overline{A_{53}}$, and thus $\Psi$ must be precisely $\IG(W_3)$, up to equivalence.
\end{proof}

\subsection{Degree $4$}
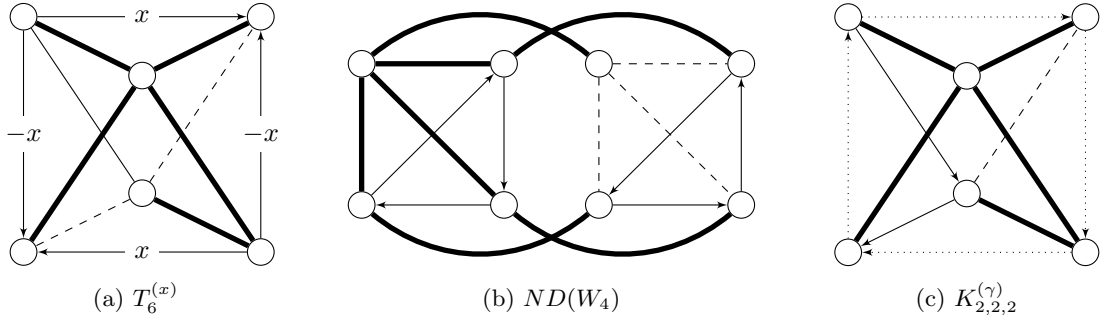
\begin{figure}[t]
    \centering
    \begin{subfigure}[b]{.25\textwidth}\centering
    \begin{tikzpicture}[scale=1.3]
        \def\x{.6}
        \def\r{2}
            \node[vertex] (1) at (0,\x) {};
            \node[vertex] (6) at (0,-\x) {};
            
            \node[vertex] (2) at (\r*\x,\r*\x) {};
            \node[vertex] (3) at (\r*\x,-\r*\x) {};
            \node[vertex] (4) at (-\r*\x,-\r*\x) {};
            \node[vertex] (5) at (-\r*\x,\r*\x) {};
            
            \draw[tree] (1) to node{} (2);
            \draw[tree] (1) to node{} (3);
            \draw[tree] (1) to node{} (4);
            \draw[tree] (1) to node{} (5);
            
            \draw[negedge] (6) to node{} (2);
            \draw[tree] (6) to node{} (3);
            \draw[negedge] (6) to node{} (4);
            \draw[edge] (6) to node{} (5);
            
            \draw[arc] (3) to node[fill=white]{$-x$} (2);
            \draw[arc] (3) to node[fill=white]{$x$} (4);
            \draw[arc] (5) to node[fill=white]{$-x$} (4);
            \draw[arc] (5) to node[fill=white]{$x$} (2);
            
    \end{tikzpicture}
    \caption{$T_6^{(x)}$}
    \label{fig: t6x}
    \end{subfigure}
    \begin{subfigure}[b]{.48\textwidth}\centering
    \begin{tikzpicture}[scale=1.2*1.3]
        \def\x{.6}
        \def\r{2}
        \def\b{40}
        
            \node[vertex] (1) at (-\x, \x) {};
            \node[vertex] (2) at (-\x, -\x) {};
            \node[vertex] (3) at (\x, \x) {};
            \node[vertex] (4) at (\x, -\x) {};
            
            \node[vertex] (5) at (\r-\x, \x) {};
            \node[vertex] (6) at (\r-\x, -\x) {};
            \node[vertex] (7) at (\r+\x, \x) {};
            \node[vertex] (8) at (\r+\x, -\x) {};
            
            \draw[tree, bend left=\b] (1) to node{} (5);
            \draw[tree, bend left=\b] (3) to node{} (7);
            \draw[tree, bend right=\b] (2) to node{} (6);
            \draw[tree, bend right=\b] (4) to node{} (8);
            
            \draw[tree] (1) to node{} (2);
            \draw[tree] (1) to node{} (3);
            \draw[tree] (1) to node{} (4);
            
            \draw[negedge] (5) to node{} (6);
            \draw[negedge] (5) to node{} (7);
            \draw[negedge] (5) to node{} (8);
            
            \draw[arc] (2) to node{} (3);
            \draw[arc] (3) to node{} (4);
            \draw[arc] (4) to node{} (2);
            
            \draw[arc] (7) to node{} (6);
            \draw[arc] (8) to node{} (7);
            \draw[arc] (6) to node{} (8);
    \end{tikzpicture}
    \caption{$ND(W_4)$}
    \label{fig: nd w4}
    \end{subfigure}
        \begin{subfigure}[b]{.25\textwidth}\centering
    \begin{tikzpicture}[scale=1.3]
        \def\x{.6}
        \def\r{2}
        
            \node[vertex] (1) at (0,\x) {};
            \node[vertex] (6) at (0,-\x) {};
            
            \node[vertex] (2) at (\r*\x,\r*\x) {};
            \node[vertex] (3) at (\r*\x,-\r*\x) {};
            \node[vertex] (4) at (-\r*\x,-\r*\x) {};
            \node[vertex] (5) at (-\r*\x,\r*\x) {};
            
            \draw[tree] (1) to node{} (2);
            \draw[tree] (1) to node{} (3);
            \draw[tree] (1) to node{} (4);
            \draw[tree] (1) to node{} (5);
            
            \draw[negedge] (6) to node{} (2);
            \draw[tree] (6) to node{} (3);
            \draw[arc] (6) to node{} (4);
            \draw[arc] (5) to node{} (6);
            
            \draw[arc, dotted] (2) to node{} (3);
            \draw[arc, dotted] (3) to node{} (4);
            \draw[arc, dotted] (4) to node{} (5);
            \draw[arc, dotted] (5) to node{} (2);

    \end{tikzpicture}
    \caption{$K_{2,2,2}^{(\gamma)}$}
    \label{fig: k222gamma}
    \end{subfigure}
    \caption{Illustrations for Proposition \ref{prop: k=4, with triangles}. Filled arcs have gain $i$. Further edges without labels follow the usual drawing conventions. }
    \label{fig: degree 4 with triangle illustrations}
\end{figure}

Proceeding along the same line, we may consider unit gain graphs with degree four. 
In order to gain some insight into the underlying graph, we use the following consequence of Lemma \ref{lemma: perfect square}.
\begin{lemma}\label{onetriangle}
Let $k$ be such that $1+4k$ is not a square, and let $\Psi$ be a two-eigenvalue gain graph with degree $k$. 
Then in the underlying graph, no edge can be in exactly one triangle.
\end{lemma}
\begin{proof} 
Assume to the contrary that $u$ and $v$ are adjacent with precisely one common neighbor $w$. 
Without loss of generality, we may assume that $A_{uv}=A_{uw}=1$. 
Then $A_{vw}=(A^2)_{uv}=a$, which is real, since $a=\theta_1+\theta_2$, and unit, since $a=A_{vw}$.
Therefore, $a=\pm1$. 
But then $a\in \mathbb{N}_0$, so  $a^2 + 4k$ is a perfect square, by Lemma \ref{lemma: perfect square}, which is a contradiction. 
\end{proof}
In the following classification, we first treat the case in which triangles are allowed. 
\begin{proposition} 
\label{prop: k=4, with triangles}
Let $D$ be a two-eigenvalue gain graph with degree $4$. 
If $D$ has triangles, then $D$ is switching isomorphic to $K_5$, $\ND(W_4)$, $K_{2,2,2}^{(\gamma)}$, or $T_{6}^{(x)}$ for some unit $x$.
\end{proposition}
\begin{proof} 
Again, the $K_5$ case is clear, so let us assume that $D$ is not switching isomorphic to a complete graph.
Then both eigenvalue multiplicities are larger than $1$, which implies that $n\geq 6$, since $n=5$ would violate the absolute bound for equiangular lines in $\mathbb{C}^2$. (See Section \ref{subsec: bounds}.)

Since there is a triangle, we obtain by Lemma \ref{onetriangle} that there must be an edge which is in precisely two triangles. 
Indeed, if every edge would be in 0 or more than 2 triangles, then the graph would be $K_5.$
Without loss of generality, we may now assume that $A_{12}=A_{13}=A_{14}=A_{15}=1$, $A_{23}=0$, $A_{24}=x\in \T$, and $A_{25}\neq 0$. Now since $A_{24}+A_{25}=(A^2)_{21}=a$, for some real $a$, $A_{25}=a-x$.
We now distinguish two cases depending on whether $a=0$ or not, in which case $a>0$. 

For the first case, assume that $a=0$ and thus $A_{25}=-x$.
Let $y,z\in\T\cup\{0\}$, and without loss of generality set $A_{35}=y$, $A_{45}=z$. 
Since $(A^2)_{31}=0$, it follows that $A_{34}=-y$. 
Moreover, since $(A^2)_{41}=0$ and $(A^2)_{51}=0$, it follows that 
\begin{equation}\label{eq: x y z equations}\overline{x}-\overline{y}+z=0 \text{~and~} -\overline{x}+\overline{y}+\overline{z}=0,\end{equation}
and hence $z+\overline{z}=0.$
This holds true in two subcases: either $z=0$ or $z=\pm i$. 

In the former subcase, $z=0$ and thus $x=y$. 
Since we may, w.l.o.g., set $A_{26}=1$, we have 
\begin{align*}
(A^2)_{32}=A_{36}-1=0~\wedge~(A^2)_{42}=A_{46}+1=0~\wedge~(A^2)_{52}=A_{56}+1=0 \\
~\implies A_{36}=-A_{46}=-A_{56}=1,
\end{align*}
and we obtain the toral tesselation graph $T_6^{(x)}$, illustrated in Figure \ref{fig: t6x}, with no further restrictions on $x$. 

In the latter subcase, w.l.o.g. choose $z=i$. 
Applying the same technique again, we find 
\[0=(A^2)_{45}=1-x\overline{x}-y\overline{y}=-y\overline{y}\implies y=0,\]
which means that $n\geq 8$, since vertex $3$ needs three more neighbors.
Moreover, by plugging in $y$ and $z$ into \eqref{eq: x y z equations}, it follows that $x=i$. 
Without loss of generality we may then assume that $A_{26}=1$ and observe that $(A^2)_{2j}=0$ for all $j$, to determine $A_{j6}$ for $j=3,4,5$. 
Similarly by assuming (w.l.o.g.) $A_{47}=1$ we may determine $A_{j7}$ for $j=3,5,6$, and finally assuming $A_{58}=1$ determines $A_{j8}$ for $j=3,6,7$.
Altogether, we find a unique graph (up to switching equivalence), which was before obtained as $\ND(W_4)$, and is illustrated in Figure \ref{fig: nd w4}.\\

In case $a>0$, we find one more switching equivalence class through a series of similar arguments.
Recall that $A_{24}=x$ and $A_{25}=a-x$.
Since $a-x\in\T$, it follows that $a-x=\overline{x}.$ 
Similarly, since $(A^2)_{13}=A_{43}+A_{53} = a = x+\overline{x},$ it follows that $A_{43}\in\{x,\overline{x}\}$ and $A_{53}=\overline{A_{43}}.$
However, setting $A_{26}\not=0$ w.l.o.g. and choosing  $A_{43}=\overline{x}$ yields 
\[(A^2)_{23}=A_{21}A_{13} + A_{24}A_{43} + A_{25}A_{53}+ A_{26}A_{63} = 3+A_{26}A_{63}\not=0=a\cdot A_{23},\]
where the inequality  holds since $|A_{26}A_{63}|\leq 1$. 
Clearly, this is a contradiction and thus $A_{43}=A_{35}=x$.  
From $(A^2)_{14}$ it then follows similarly that $A_{45}=0$.

As before, we may now assume without loss of generality that $A_{26}=1$, and determine the values $A_{63}, A_{64}$ and $A_{65}$ by repeating the same argument three times, as follows. 
\begin{align*}
    1 + A_{64} = (A^2)_{24} = a\cdot A_{24} = (x+\overline{x}) x = 1 + x^2 &\implies A_{64}=x^2, \\
    1 + A_{65} = (A^2)_{25} = a\cdot A_{25} = (x+\overline{x}) \overline{x} = 1 + \overline{x}^2 &\implies A_{65}=\overline{x}^2, \text{~and} \\
    1 + \overline{x}^2A_{63} = (A^2)_{43} = a\cdot A_{43} = (x+\overline{x}) x = 1 + x^2 &\implies A_{63}=x^4.
\end{align*}
Finally, note that  $(A^2)_{61} = 1+x^2+\overline{x}^2+x^4=(1+\overline{x}^2)(1+x^4)=0$, which holds subject to $a=x+\overline{x}>0$ precisely when either $x=\gamma$ or $x=\overline{\gamma}$; both cases yield $K_{2,2,2}^{(\gamma)}$, illustrated above. 
\end{proof}

The attentive reader may have noted that we did not fix the order of the above considered graphs, during the proof. 
However, by allowing the initially undetermined edge gains to be either complex units or zero, and following the implications, we arrive at a $4$-regular gain graph in each of the possible cases. 
Hence, we may be certain that no larger connected examples exist. 

What is left is the triangle-free case. 
While the following result may be shown through procedure similar to the above, it has been shown by \cite[Sec. 3.5]{best2013unit}, in the classification of weighing matrices.
Indeed, note that by Lemma \ref{lemma: triangle-free implies a=0}, the gain matrices of the examples below are all weighing matrices. 
We thus omit the rather lengthy proof. 
\begin{proposition}
\label{prop: k=4, no triangle}
Let $D$ be a two-eigenvalue gain graph with degree $4$. 
If $D$ is triangle-free, then $D$ is switching isomorphic to one of the graphs $\IG(W_5)$, $\ND(\IG(W_3))$, $\IG(W_7)$, $\ND(\ND(\IG(W_2)))$, or $T_{2t}^{(x)}$, for some $t>3$ and $x\in \T$.
\end{proposition}

The results of our classification are summarized by the following theorem. 
\begin{theorem}
\label{thm: summary}
All two-eigenvalue gain graphs with degree at most $4$ are switching isomorphic to one of the gain graphs in Table \ref{tab: summary}.
\end{theorem}
\begin{table}[h!]
\centering
\begin{tabular}{@{}clllccclllc@{}}
\toprule
$k$ & \multicolumn{1}{l}{Graph} & \multicolumn{1}{l}{Order} & \multicolumn{1}{l}{$m$} &\multicolumn{1}{l}{DS} && $k$ & \multicolumn{1}{l}{Graph} & \multicolumn{1}{l}{Order} & \multicolumn{1}{l}{$m$}&\multicolumn{1}{l}{DS} \\ \cmidrule{1-5} \cmidrule{7-11}
\multirow{2}{*}{$2$}  & $K_3$                    & $3$                       & $1$      &  *          &&   
 \multirow{8}{*}{$4$} & $K_5$                     & $5$                      & $1$      &   *               \\
                     & $IG(W_2)$                 & $4$                       & $2$      &  *          &&     
                      & $K_{2,2,2}^{(\gamma)}$    & $6$                        & $2$    & *                 \\
                     &                           &                           &          &            &&     
                      & $ND(W_4)$                 & $8$                       & $4$   &  *                \\ 
                      &                           &                           &          &            &&     
                      & $IG(W_5)$                 & $10$                       & $5$   &  *                \\
 \multirow{4}{*}{$3$} & $K_4$                     & $4$                       & $1$     &  *          &&     
                      & $ND(IG(W_3))$             & $12$                      & $6$  &  *                \\ 
                      & $W_4$                     & $4$                       & $2$  &  *          &&     
                      & $IG(W_7)$                 & $14$                      & $7$  &  *                \\ 
                      & $IG(W_3)$                 & $6$                       & $3$  &  *          &&     
                      & $ND(ND(IG(W_2)))$         & $16$                      & $8$  &  *                \\ 
                      & $ND(IG(W_2))$             & $8$                       & $4$  &  *          &&     
                      & $T_{2t}^{(x)}$            & $2t,~t\geq 4$             & $t$      \\ \bottomrule
\end{tabular}
\caption{Classification of all two-eigenvalue gain graphs with degree at most $4$. A star in the DS column indicates that any connected, cospectral gain graph is switching isomorphic.}
\label{tab: summary}
\end{table}
\subsection{Degree $5$}
\label{sec: degree 5}
To conclude the discussion of low-degree two-eigenvalue gain graphs, we consider some degree-$5$ examples. 
We briefly touch on graphs of order at most $8$, after which we will treat a new infinite family of two-eigenvalue gain graphs, whose underlying structure is somewhat of a doubled cycle. 
The remaining sporadic examples that have been found through computer search appear in Appendix \ref{sec: SA results}, and will not be discussed explicitly. 

Since there are exactly four $5$-regular graphs of order at most $8$, we may simply treat them on a case-by-case basis. 
For two of those candidates, namely the complement of $C_8$, the complement of $C_3\cup C_5$, one may show that neither may be underlying to a two-eigenvalue gain graph. 
The proof follows the same pattern as the proofs of Propositions \ref{prop: k=4, with triangles} and \ref{prop: k=4, no triangle}, so we forego the details.   

\begin{figure}[t]
    
    \centering
    \begin{subfigure}[b]{.48\textwidth}\centering
    \begin{tikzpicture}
        \def\x{.6}
        \def\r{2}
        \def\n{6}
        \def\lw{1}
        \foreach \a in {1,2,...,\n}{
            \node[vertex] (o\a) at (180+\a*360/7: 2 cm) {};
            }
            \foreach \a in {1,2,...,\n}{
            \node[vertex] (i\a) at (180+\a*360/7: 1 cm) {};
            }
            \node[vertex] (i0) at (-1,0) {};
            \node[vertex] (o0) at (-2,0) {};
            
            \draw[edge,dotted] (i0) to node{} (i1);
            \draw[edge] (i1) to node{} (i2);
            \draw[edge] (i2) to node{} (i3);
            \draw[arc,line width=\lw] (i3) to node{} (i4);
            \draw[edge] (i4) to node{} (i5);
            \draw[edge] (i5) to node{} (i6);
            \draw[edge] (i6) to node{} (i0);
            
            \draw[edge,dotted] (o0) to node{} (o1);
            \draw[negedge] (o1) to node{} (o2);
            \draw[negedge] (o2) to node{} (o3);
            \draw[negarc, line width=\lw] (o3) to node{} (o4);
            \draw[negedge] (o4) to node{} (o5);
            \draw[negedge] (o5) to node{} (o6);
            \draw[negedge] (o6) to node{} (o0);
            
            \draw[edge] (i0) to node{} (o0);
            \draw[edge] (i1) to node{} (o1);
            \draw[edge] (i2) to node{} (o2);
            \draw[edge] (i3) to node{} (o3);
            \draw[edge] (i4) to node{} (o4);
            \draw[edge] (i5) to node{} (o5);
            \draw[edge] (i6) to node{} (o6);

            \draw[edge] (i1) to node{} (o2);
            \draw[edge] (i2) to node{} (o3);
            \draw[arc,line width=\lw] (i3) to node{} (o4);
            \draw[edge] (i4) to node{} (o5);
            \draw[edge] (i5) to node{} (o6);
            
            \draw[negedge] (o1) to node{} (i2);
            \draw[negedge] (o2) to node{} (i3);
            \draw[negarc,line width=\lw] (o3) to node{} (i4);
            \draw[negedge] (o4) to node{} (i5);
            \draw[negedge](o5) to node{} (i6);

            \draw[edge,dotted] (i0) to node{} (o1);
            \draw[negedge] (i0) to node{} (o6);
            
            \draw[edge,dotted] (o0) to node{} (i1);
            \draw[edge] (o0) to node{} (i6);
    \end{tikzpicture}
    \caption{A general two-eigenvalue donut}
    \label{fig: general donut}
    \end{subfigure}
    \begin{subfigure}[b]{.48\textwidth}\centering
    \begin{tikzpicture}[scale=1.25]
        \def\x{.5}
        \def\r{3.1}
        
            \node[vertex] (1) at (\x,\x) {};
            \node[vertex] (2) at (\x,-\x) {};
            \node[vertex] (3) at (-\x,-\x) {};
            \node[vertex] (4) at (-\x,\x) {};
            
            \node[vertex] (5) at (\r*\x,\r*\x) {};
            \node[vertex] (6) at (\r*\x,-\r*\x) {};
            \node[vertex] (7) at (-\r*\x,-\r*\x) {};
            \node[vertex] (8) at (-\r*\x,\r*\x) {};
            
            \draw[edge] (1) to node{} (2);
            \draw[edge] (3) to node{} (2);
            \draw[edge] (1) to node{} (4);
            \draw[edge] (3) to node{} (4);
            
            \draw[negedge] (5) to node{} (6);
            \draw[negedge] (7) to node{} (6);
            \draw[negedge] (5) to node{} (8);
            \draw[negedge] (7) to node{} (8);
            
            \draw[edge] (1) to node{} (5);
            \draw[edge] (2) to node{} (6);
            \draw[edge] (3) to node{} (7);
            \draw[edge] (4) to node{} (8);
            
            \draw[arc] (1) to node[fill=white]{$c$} (6);
            \draw[arc] (7) to node[fill=white]{$c$} (2);
            \draw[arc] (3) to node[fill=white]{$c$} (8);
            \draw[arc] (5) to node[fill=white]{$c$} (4);
            
            \draw[arc] (1) to node[fill=white]{$-c$} (8);
            \draw[arc] (5) to node[fill=white]{$-c$} (2);
            \draw[arc] (3) to node[fill=white]{$-c$} (6);
            \draw[arc] (7) to node[fill=white]{$-c$} (4);

    \end{tikzpicture}
    \caption{$D^*_8(c)$}
    \label{fig: dstar8}
    \end{subfigure}
    \caption{Illustrations for Theorem \ref{thm: donuts}. The fat arcs have gain $\pm x$, and dotted lines indicate continuation of the pattern.}
    \label{fig:donuts}
\end{figure}
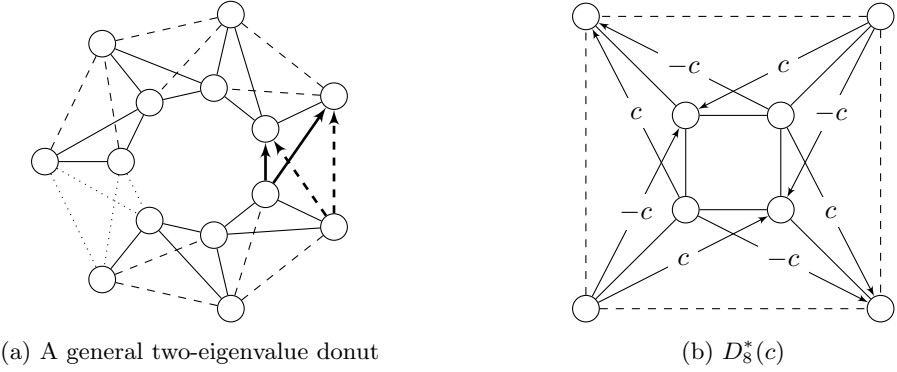

The two remaining candidates, $K_6$ and the complement of $2C_4$, are contained in the class of donut graphs, which admit infinitely many two-eigenvalue gain graphs for every (even) order $n$. 
Let us provide the formal definition. 
\begin{definition}\label{def: donut}
Let $C_t$ be the cycle graph of order $t\geq 3$, whose adjacency matrix is $B$. Then the $5$-regular graph characterized by 
\[A = \begin{bmatrix} B & B+I\\ B+I & B\end{bmatrix}\]
is called the order-$2t$ \textit{donut graph.}
\end{definition}
We will now characterize all two-eigenvalue donut graphs with symmetric spectra. 
\begin{theorem}\label{thm: donuts}
Let $G$ be an order $n:=2t$ donut graph, $t\geq 3$, and let $\Psi = (G,\psi)$ be a unit gain graph. 
Then $\Psi$ has eigenvalues $\pm \sqrt{5}$ if and only if it is switching isomorphic to 
\[\begin{bmatrix} C + C^* & C-C^* + I \\ C^*-C+I & -C-C^*\end{bmatrix},\]
where $C$ is an order-$t$ weighing matrix of weight $1$, or to $D^*_8(c)$ in Figure \ref{fig: dstar8}.
\end{theorem}

\begin{proof}
Sufficiency is clear, following the discussion in Section \ref{sec: constructions}, so we only show necessity. 
If $t=3$ then the claim holds by Theorem \ref{thm: mult at most 3}, so suppose that $t\geq 4$. 
Let $C$ be the $t\times t$ matrix whose nonzero entries are $C_{jh}=1$ for all $h=j+1$ and $C_{t1}=x$ for some $x\in T$. 
Let $W,Y$ and $Z$ be matrices with the same support as $C$, and set 
\[A = \begin{bmatrix}C+C^* & W-Z^*+I\\ W^*-Z+I & Y+Y^*\end{bmatrix}.\]
Note that without loss of generality, $A(\Psi)=A$. 
Computing the upper left block of $A^2(\Psi)$ yields:
\begin{align*}
A^2=5I &~\implies ~ C^2 - WZ + (C^*)^2  - Z^*W^* + W - Z  + W^*- Z^*  = O \\
 &~\implies ~  W=Z \wedge C^2 + (C^*)^2= W^2 + (W^*)^2, 
\end{align*}
where the final equivalence follows by grouping the terms by support. 
Similarly, plugging in the above and computing the upper right block yields 
\begin{align*}
 A^2=5I 
 &~\implies~ C + Y + C^* + Y^* + CW + WY -C^*W^*- W^*Y^*  = O \\
 &~\implies~  Y=-C \wedge CW + W^*C^*=WC + C^*W^*.
\end{align*}
Now, if $t\geq 5$ then we may again group by support to reduce the final equality above to $CW=WC$.
It follows that either $W=C$ or $W=-C$, completing the gain graph. 
Note that either choice yields the same switching equivalence class.

In case $t=4$, all of the second order matrices have the same supports, so the above is not the only solution. 
Briefly put, by solving the system 
\[\begin{cases}C^2 + (C^*)^2= W^2 + (W^*)^2 \\ CW + W^*C^*=WC + C^*W^*
\end{cases}\iff
\begin{cases}
\overline{x} + 1 = w_{12}w_{23} + \overline{w_{34}}\overline{w_{41}}\\
\overline{x} + 1 = w_{23}w_{34} + \overline{w_{12}}\overline{w_{41}}\\
w_{23} + \overline{w_{41}} = w_{12} + \overline{w_{34}}\overline{x}\\
w_{34} + \overline{w_{12}}\overline{x} = w_{23} + \overline{w_{41}}
\end{cases}\]
we obtain either $W=\pm C$, as above, or either one of
\[x=-1 ~\wedge W = \text{diag}\left(\begin{bmatrix}c& \bar{c} & \bar{c} & c\end{bmatrix}\right)C \text{~~and~~}
 x=1 ~\wedge W = \text{diag}\left(\begin{bmatrix}c& \bar{c} &c & \bar{c}\end{bmatrix}\right)C,\]
where $c\in\T$.
The latter two yield gain graphs switching isomorphic to the exception $D_8^*$, shown in Figure \ref{fig: dstar8}, which completes the proof. 
\end{proof}
Note that indeed, $K_6$ is a donut, stringly speaking; though it is somewhat of a special case. 
In particular, since $t=3$ implies that $C^2$ and $C^*$ have the same support, the particulars of the proof above do not apply. 
However, as a consequence of Theorem \ref{thm: mult at most 3}, we know that the statement holds regardless.

Finally, note that the case $D_8^*(c)$ is distinct from the two-eigenvalue order-$8$ donut that follows the general construction. 
Indeed, the triangles in the former have gains $\pm c$, whereas the latter has triangles with gains $\pm 1$ for all $x\in\T$.

\section{Computer-aided search}
\label{sec: heuristic}

In addition to the various sources in the literature that are related to two-eigenvalue gain graphs, it has proved fruitful to develop a computer-aided search method.
In more traditional graph related fields, one is often inclined to iteratively consider all graphs within a given band of parameters. 
However, since we here consider essentially weighted graphs with weights on the complex unit circle, this approach does not work for the current context. 
Indeed, note that we may not simply restrict to (sparse) multiplicative subsets of $\T$, without potentially missing examples, such as those from Theorem \ref{thm: renes family}. 

Nonetheless, we have been able to implement the somewhat unusual, though generally well-known Simulated Annealing \cite{kirkpatrick1983optimization} procedure to fit our needs. 
Let $f(\cdot)$ be a function that measures the quality of the proposed candidates. 
For this particular application, we found the the function 
\[f(A) = \|A^2 - (\lambda_1+\lambda_n)A + \lambda_1\lambda_nI\|,\]
where $\lambda_1$ and $\lambda_n$ are respectively the largest and the smallest eigenvalues of $A$, suited our needs well.
Broadly speaking, the procedure for a given graph may may then be described as follows. 
Here, $\alpha<1$ is the cool-down parameter, and $\tau$ is the minimum temperature.  
\begin{enumerate}
    \item[Init)]    Fix a spanning tree of the edges, and randomly assign the other edges a gain.
    \item[(A)]        Perform the following $m$ times:
    \begin{enumerate}
        \item[i.]    Randomly rotate every non-fixed gain in $old$ at most $180t$ degrees in either direction along the complex unit circle to obtain $new$.
        \item[ii.]    Compute $f(old)$ and $f(new)$. If $f(new)<\varepsilon$, accept the candidate and stop. Else,
          set $old\leftarrow new$ with probability $\min\left\{\exp\left((f(old)-f(new))/(f(old)\cdot t)\right),1\right\}$.
    \end{enumerate}
          \item[(B)] Set $t\leftarrow \alpha t.$ If $t>\tau,$ go to (A). Otherwise, stop.
\end{enumerate}

Even when the value of the function $f$ is reasonably low, at the end of the procedure above, it yields gain graphs whose eigenvalues \textit{approximate} one of two values; one still has to distill the actual exact gain graph, to which the procedure had been converging.

We explicitly note that the above is, in essence, a clever random search algorithm.
Hence, one cannot be sure that failure to produce a two-eigenvalue gain graph $\Psi$ on a given input graph $G$ means that there are no $\Psi$ to be found on $G$.
However, some very strong theoretical convergence results have been shown \cite{henderson2003theory}. 
We have tested the success rate of this procedure on $\Gamma(\Psi)$ for various known two-eigenvalue gain graphs $\Psi$ of increasing order and density.
The results are shown below, in Figure \ref{fig: SA success}.

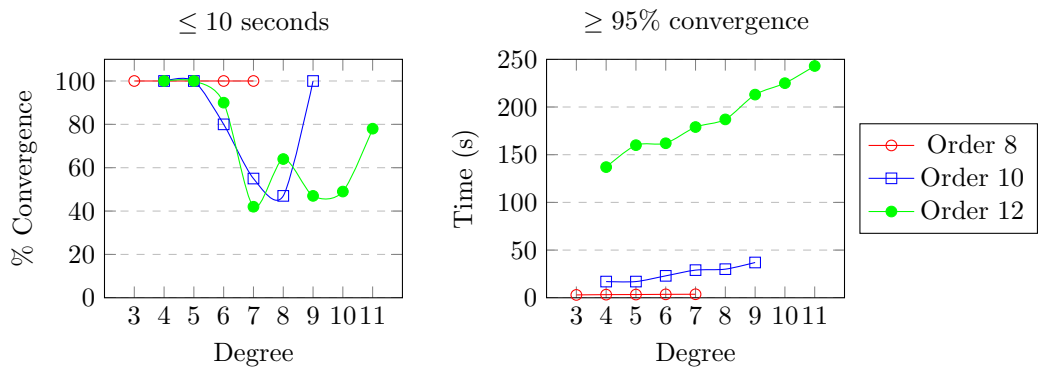
\begin{figure}[b]
    \centering
      \begin{subfigure}[b]{.40\textwidth}\centering
\pgfplotsset{width=5.5cm,compat=1.9}
\begin{tikzpicture}
\begin{axis}[
    title={$\leq 10$ seconds},
    xlabel={Degree},
    ylabel={\% Convergence},
    xmin=2, xmax=12,
    ymin=0, ymax=110,
    xtick={3,4,5,6,7,8,9,10,11,100},
    ytick={0,20,40,60,80,100},
	legend style={at={(1.05,0.5)},
	anchor=west},
    ymajorgrids=true,
    grid style=dashed,
]

\addplot[
    smooth,
    color=red,
    mark=o,
    ]
	coordinates {
	(3,100) 
	(4,100)
	(5,100) 
	(6,100) 
	(7,100)};

\addplot[
    smooth,
    color=blue,
    mark=square,
    ]
coordinates {
	(4,100) 
	(5,100) 
	(6,80) 
	(7,55) 
	(8,47) 
	(9,100)};
    
\addplot[
    smooth,
    color=green,
    mark=*,
    ]
	coordinates {
	(4,100) 		
	(5,100) 
	(6,90) 
	(7,42) 
	(8,64) 
	(9,47) 
	(10,49) 
	(11,78)};
    
\end{axis}
\end{tikzpicture}
    \end{subfigure}
      \begin{subfigure}[b]{.55\textwidth}\centering
      \pgfplotsset{width=5.5cm,compat=1.9}
\begin{tikzpicture}
\begin{axis}[
    title={$\geq 95\%$ convergence},
    xlabel={Degree},
    ylabel={Time (s)},
    xmin=2, xmax=12,
    ymin=0, ymax=250,
    xtick={3,4,5,6,7,8,9,10,11,100},
    ytick={0,50,100,150,200,250},
	legend style={at={(1.05,0.5)},
	anchor=west},
    ymajorgrids=true,
    grid style=dashed,
]

\addplot[
    smooth,
    color=red,
    mark=o,
    ]
	coordinates {
	(3,3) 
	(4,3.3)
	(5,3.4) 
	(6,3.6) 
	(7,3.7)};

\addplot[
    smooth,
    color=blue,
    mark=square,
    ]
coordinates {
	(4,17) 
	(5,17) 
	(6,23) 
	(7,29) 
	(8,30) 
	(9,37)};
    
\addplot[
    smooth,
    color=green,
    mark=*,
    ]
	coordinates {
	(4,137) 		
	(5,160) 
	(6,162) 
	(7,179) 
	(8,187) 
	(9,213) 
	(10,225) 
	(11,243)};
    \legend{Order 8, Order 10, Order 12}    
    
\end{axis}
\end{tikzpicture}
      \end{subfigure}
    
    \caption{Statistics on simulated annealing calls on admissible underlying graphs that converge to a gain graph with two eigenvalues. (40 calls for each admissible graph.)}
    \label{fig: SA success}
\end{figure}

The above shows two distinct set-ups of the various tuning parameters, one has available when applying simulated annealing. 
We primarily tune with \textit{initial temperature}, \textit{iterations per temperature}  and the \textit{cool-down procedure}, though they all more or less accomplish the same thing: more iterations of the random search, at some point in the temperature curve.

The leftmost set-up is designed to be relatively quick, with a maximum running time of at most 10 seconds. 
Conversely, the rightmost set-up is intended to yield a very high accuracy, at the expense of time. 
With the exception of the complete graphs, that somehow converge remarkably well in the quick set-up, we observe a clear trade-off.

It should also be noted that almost exactly the same routine can be used to search for cospectral mates of arbitrary gain graphs. 
Indeed, all one really has to change is the evaluation function $f$, which may simply be set to the sum of the quadratic differences with respect to some predetermined spectrum. 
For sparse graphs in particular, this idea could be easily applied to get a feel for the distinct, cospectral gain graphs that satisfy a given mix of structural properties, if one were to investigate, say, spectrally determined gain graphs. 


\scriptsize
\DeclareRobustCommand{\VAN}[3]{#3}
\bibliographystyle{abbrvnat}
\bibliography{mybib}{}
\normalsize
\appendix
\section{Simulated Annealing Results}\label{sec: SA results}
This section contains some of the results that have been found by implementation of the algorithm discussed in Section \ref{sec: heuristic}, which have not appeared elsewhere in this work. 
In a systematic search of two-eigenvalue gain graphs on order-$12$, degree-$5$ graphs, we found precisely $M_1,\ldots,M_4$, below, in addition to the order-$12$ donut. 
As mentioned before, we cannot be certain that no other examples exist, due to the random nature of the search algorithm, but the authors are reasonably confident that all two-eigenvalue gain graphs with $(n,k)=(12,5)$ have been found. 
\begin{itemize} 
\item  New example on $K_8$:
\begin{equation*}K_8^*=\begin{bmatrix}
 0&   1&   1&   1&   1&   1&   1&   1\\
 1&   0&  i& -i&  i& -i&  i& -i\\
 1& -i&   0& -i& -i&  i&  i&  i\\
 1&  i&  i&   0& -i& -i& -i&  i\\
 1& -i&  i&  i&   0&  i& -i& -i\\
 1&  i& -i&  i& -i&   0&  i& -i\\
 1& -i& -i&  i&  i& -i&   0&  i\\
 1&  i& -i& -i&  i&  i& -i&   0
\end{bmatrix}.\end{equation*}

\item A signed graph example on $K_{10}$:
\begin{equation*}K_{10}^*=\begin{bmatrix}
      0 &    1 &    1 &    1 &    1 &    1 &    1 &    1 &    1 &    1 \\
      1 &    0 &   -1 &   -1 &    1 &    1 &   -1 &   -1 &    1 &    1\\
      1 &   -1 &    0 &    1 &    1 &   -1 &    1 &   -1 &   -1 &    1\\
      1 &   -1 &    1 &    0 &   -1 &   -1 &   -1 &    1 &    1 &    1\\
      1 &    1 &    1 &   -1 &    0 &   -1 &    1 &   -1 &    1 &   -1\\
      1 &    1 &   -1 &   -1 &   -1 &    0 &    1 &    1 &   -1 &    1\\
      1 &   -1 &    1 &   -1 &    1 &    1 &    0 &    1 &    -1 &   -1\\
      1 &   -1 &   -1 &    1 &   -1 &    1 &    1 &    0 &    1 &   -1\\
      1 &    1 &   -1 &    1 &    1 &   -1 &   -1 &    1 &    0 &   -1\\
      1 &    1 &    1 &    1 &    -1 &    1 &   -1&    -1&    -1 &   0
\end{bmatrix}.\end{equation*}

\newpage
\item New example on the icosahedron:

\setcounter{MaxMatrixCols}{20}
\begin{equation*}M_1 = \begin{bmatrix} 0&  1&  0&  0&  0&  1&  1&  1&  0&  0&  0&  1\\
 1&  0&  1&  0&  0&  0&  1& -1&  1&  0&  0&  0\\
 0&  1&  0&  1&  0&  0&  0& -1& -1&  1&  0&  0\\
 0&  0&  1&  0&  1&  0&  0&  0& -1& -1&  1&  0\\
 0&  0&  0&  1&  0&  x&  0&  0&  0& -1& -1& -x\\
 1&  0&  0&  0& \bar{x}&  0& -1&  0&  0&  0&-\bar{x}&  1\\
 1&  1&  0&  0&  0& -1&  0&  1&  0&  0&  0& -1\\
 1& -1& -1&  0&  0&  0&  1&  0& -1&  0&  0&  0\\
 0&  1& -1& -1&  0&  0&  0& -1&  0& -1&  0&  0\\
 0&  0&  1& -1& -1&  0&  0&  0& -1&  0& -1&  0\\
 0&  0&  0&  1& -1& -x&  0&  0&  0& -1&  0&  x\\
 1&  0&  0&  0&-\bar{x}&  1& -1&  0&  0&  0& \bar{x}&  0 \end{bmatrix},~x\in\T.\end{equation*}

\item A bipartite example based on a novel non-graphical weighing matrix of weight 5: 
\begin{equation*}
Z=\begin{bmatrix}
1&    1&     1&     1&      1&    0 \\ 
1&   -\bar{x}&    -1&     \bar{x}&     0&   -x  \\
1&    -1&     x&     0&     -x&    x  \\
1&    \bar{x}&     0&    -\bar{x}&    -1&   -x  \\
1&     0&    -x&    -1&      x&    x  \\
0&    \bar{x}&    -\bar{x}&    \bar{x}&   -\bar{x}&    1  
\end{bmatrix}, x\in\T,  \text{ and } M_2 = \begin{bmatrix} O & Z \\ Z^* & O\end{bmatrix}.
\end{equation*}

\item Two more sporadic examples:
\begin{equation*}
    M_3=\begin{bmatrix}
  0&  0&        0&        0&        1&  0&  0&        1&  1&  0&  1&  1\\
  0&  0&        0&        0&        1&  0&  0&        1&  0&  1& -1& -1\\
  0&  0&        0&        0&        0&  1&  1&        0& -x&  x&  x&  0\\
  0&  0&        0&        0&        0&  1&  1&        0&  x& -x&  0& -x\\
  1&  1&        0&        0&        0&  0&  1&        0&  0&  0& -x&  x\\
  0&  0&        1&        1&        0&  0&  0&       -1&  1&  1&  0&  0\\
  0&  0&        1&        1&        1&  0&  0&        0& -1& -1&  0&  0\\
  1&  1&        0&        0&        0& -1&  0&        0&  0&  0&  x& -x\\
  1&  0&      -\bar{x}&       \bar{x}&        0&  1& -1&        0&  0&  0&  0&  0\\
  0&  1&       \bar{x}&      -\bar{x}&        0&  1& -1&        0&  0&  0&  0&  0\\
  1& -1&       \bar{x}&        0&      -\bar{x}&  0&  0&       \bar{x}&  0&  0&  0&  0\\
  1& -1&        0&      -\bar{x}&       \bar{x}&  0&  0&      -\bar{x}&  0&  0&  0&  0\\
    \end{bmatrix}, x\in \T, \text{ and }
\end{equation*}
\begin{equation*}
    M_4 = \begin{bmatrix}       0 &    0 &    0 &    1 &    0 &    0 &    1 &    0 &    1 &    1 &    1 &    0 \\
   0 &    0 &    0 &    0 &    1 &    0 &    1 &    1 &    0 &    0 &   -1 &    1\\ 
   0 &    0 &    0 &    0 &    0 &    1 &    1 &   -1 &    0 &    0 &   -1 &   -1 \\ 
   1 &    0 &    0 &    0 &    0 &    0 &    i &    0 &   -i &   -i &    i &    0 \\ 
   0 &    1 &    0 &    0 &    0 &    0 &    0 &   -i &    i &   -i &    0 &    i \\ 
   0 &    0 &    1 &    0 &    0 &    0 &    0 &    i &    i &   -i &    0 &   -i \\ 
   1 &    1 &    1 &    -i &    0 &    0 &    0 &    0 &    0 &    i &    0 &    0 \\ 
   0 &    1 &   -1 &    0 &   i &    -i &    0 &    0 &    0 &    0 &    0 &   -i \\ 
   1 &    0 &    0 &   i &    -i &    -i &    0 &    0 &    0 &    0 &   -i &    0 \\ 
   1 &    0 &    0 &   i &   i &   i &    -i &    0 &    0 &    0 &    0 &    0 \\
   1 &   -1 &   -1 &    -i &    0 &    0 &    0 &    0 &   i &    0 &    0 &    0 \\ 
   0 &    1 &   -1 &    0 &    -i &   i &    0 &   i &    0 &    0 &    0 &    0 
   \end{bmatrix}.
\end{equation*}
\end{itemize}

\end{document}